\documentclass[11pt]{article}
\usepackage{a4}
\usepackage{graphicx}
\usepackage{parskip}

\usepackage{amsfonts}
\usepackage{natbib}
\usepackage{algorithm}
\usepackage{algorithmic}

\usepackage{latexsym}

\newtheorem{theorem}{Theorem}[section]

\newtheorem{lemma}{Lemma}[section]
\newtheorem{corollary}{Corollary}[section]
\newtheorem{remark}{Remark}[section]
\newtheorem{example}{Example}[section]

\newcommand{\R}{\mathbb{R}}
\newcommand{\Nat}{\mathbb{N}}

\newcommand{\PP} {{  \rm I\hskip-0.22em P}}
\newcommand{\bx}{{\bf X}}

\newcommand{\by}{{\bf Y}}

\newcommand{\EE} {{\rm I\hskip-0.48em E}}

\usepackage{natbib}

\begin{document}
\centerline {\bf \Large The adaptive and the thresholded Lasso}

\centerline {\bf \Large 
for potentially misspecified models}

\centerline{July, 2010}

\centerline{Sara van de Geer*, Peter B\"uhlmann*, and Shuheng Zhou**}

\centerline{*Seminar for Statistics, ETH Z\"urich}

\centerline{**Department of Statistics, University of Michigan}

{\bf Abstract} We revisit the adaptive Lasso as well as the thresholded Lasso with refitting, in a high-dimensional linear model, and
study  prediction error, $\ell_q$-error ($q \in \{ 1 , 2 \} $), and number of false positive selections.
Our theoretical results for the two methods
are, at a rather fine scale, comparable.
The differences only
show up in terms of the (minimal) restricted and sparse eigenvalues, favoring thresholding
over the adaptive Lasso.
As regards prediction and estimation, the difference is virtually negligible, 
but our bound for the number of false positives is larger for the
adaptive Lasso than for thresholding.
Moreover, both these 
two-stage methods add value to  the one-stage Lasso in the sense
that, under appropriate restricted and sparse eigenvalue conditions,
they have similar prediction and estimation error as the one-stage
Lasso, but substantially less false positives.


{\bf Running Head:} {Adaptive and thresholded Lasso}

\section{Introduction}
\label{introduction.section}

 Consider the linear model 
$$
{\bf Y} = {\bf X} \beta + \epsilon,
$$
where $\beta \in \R^p$ is a vector of coefficients,
${\bf X}$ is an $(n \times p)$-design matrix, and ${\bf Y}$ is an $n $-vector
of noisy observations, $\epsilon$ being the noise term. 
We examine the case $p \ge n$, i.e., a high-dimensional situation.
The design
matrix ${\bf X}$ is treated as fixed, and 
the Gram matrix is denoted by $\hat \Sigma := {\bf X}^T {\bf X}/n $. 
Throughout, we assume the
normalization
$ \hat \Sigma_{j,j} = 1$ for all $j \in \{ 1 , \ldots , p \}$.

This paper presents a theoretical comparison between the thresholded
Lasso with refitting and the adaptive Lasso. Both methods are very
popular in practical applications for reducing the number of
active variables.

We emphasize here and describe later that we allow for model
misspecification where the true regression function may be non-linear in
the covariates. For such cases, we can consider the projection onto the
linear span of the covariates. The (projected or true) linear model does
not need to be sparse nor do we require that the non-zero regression
coefficients (from a sparse approximation) are ``sufficiently large''.  As
for the latter, we will show in Lemma \ref{min.lemma} how this can be invoked to improve
the result. Furthermore, we also do not require the stringent 
irrepresentable conditions or incoherence assumptions on the 
design matrix ${\bf X}$ but only some weaker restricted or sparse
eigenvalue conditions.

Regularized estimation with the $\ell_1$-norm penalty, also 
known as the Lasso \break (\cite{Tib96}), refers to the following convex
optimization problem:   
\begin{eqnarray}
\label{eq.init}
\hat \beta := \arg\min_{\beta}\biggl \{  \|{\bf Y} -{\bf X} \beta\|_2^2/n +  \lambda \|\beta\|_1
\biggr \} ,
\end{eqnarray}
where $\lambda >0
$ is a penalization parameter. 

Regularization with $\ell_1$-penalization in high-dimensional scenarios has become extremely
popular. The methods are
easy to use, due to recent progress in specifically tailored convex optimization
(\cite{mevdgpb08}, \cite{fht08}).

A two-stage version of the Lasso is the so-called adaptive Lasso 
\begin{eqnarray}
\label{eq.adap}
\hat \beta_{\rm adap} := \arg\min_{\beta} 
\biggl \{ \|{\bf Y} -{\bf X} \beta\|_2^2/n +  \lambda_{\rm init} \lambda_{\rm adap} 
\sum_{j=1}^p {| \beta_j | \over | \hat \beta_{j, {\rm init} } |} \biggr \}  .
\end{eqnarray}
Here, $\hat {\beta}_{\rm init}$ is the one-stage Lasso defined in
(\ref{eq.init}), with initial tuning parameter $\lambda= \lambda_{\rm init}$, and
$\lambda_{\rm adap}>0$ is the tuning parameter for the second stage.
Note that when $|\hat \beta_{j , \rm init}|=0$, we exclude variable $j$ in the second stage.
The adaptive Lasso was originally
proposed by \cite{Zou06}. 

Another possibility is the thresholded Lasso with refitting. Define
\begin{eqnarray} 
\hat{S}_{\rm thres} = \{ j : |\hat \beta_{j, {\rm init} } |  > \lambda_{\rm
  thres}\}, 
\end{eqnarray}
which is the set of variables having estimated coefficients larger than
some given threshold $\lambda_{\rm thres}$. The refitting is then done by ordinary least
squares:
\begin{eqnarray*}
\hat{b}_{\rm thres} = \arg \min _{\beta_{\hat{S}_{\rm thres}}}
\|\by - \bx \beta_{\hat{S}_{\rm thres}}\|_2^2/n,
\end{eqnarray*}
where, for a set $S \subset \{ 1 , \ldots , p \}$, $\beta_{S}$ has coefficients different from zero at the components in
$S$ only.  

We will present bounds for the prediction error, its $\ell_q$-error ($q \in \{ 1 ,2 \} $), and
the number of false positives. The bounds for the two methods are qualitatively the same.
A  difference is that our variable selection
properties results for  the adaptive Lasso depend on its prediction error,
whereas for the thresholded Lasso, variable selection can be studied without
reference to its prediction error. In our analysis this leads
to a bound for the number of false positives of
the thresholded Lasso that is smaller than the one for the adaptive Lasso, 
when restricted
or sparse minimal 
eigenvalues are small and/or sparse maximal eigenvalues are large. 

Of course, such comparisons depend on how the
tuning parameters are chosen. Choosing these by cross validation is in our view
the most appropriate, but it is beyond the scope of this 
paper to present a mathematically rigorous
theory for the cross validation scheme for the adaptive and/or  thresholded Lasso
(see \cite{arlot2010survey} for a recent survey on cross validation).

\subsection{Related work}

Consistency results for the prediction error of the Lasso can be found in \cite{GR04}.
The prediction error 
 is asymptotically oracle optimal under certain conditions on the design matrix
  ${\bf X}$,
 see e.g.\ \cite{Bunea:06, Bunea:07a, Bunea:07b}, \cite{vandeG08}, \cite{bickel2009sal},
\cite{koltch09a, koltch09b},
where also estimation in terms of the $\ell_1$- or $\ell_2$-loss is considered.
The ``restricted eigenvalue condition" of \cite{bickel2009sal} (see also
\cite{koltch09a, koltch09b}) plays a key role here. 
Restricted eigenvalue conditions are implied by, but generally much weaker than,  ``incoherence" conditions,
which exclude high correlations between co-variables.
Also \cite{CP09} allow for a major relaxation of incoherence conditions, using assumptions on the set of true coefficients.

There is however a bias problem with $\ell_1$-penalization, due to the shrinking of
the estimates which correspond to true signal variables.
A discussion  
can be found in \cite{Zou06}, and \cite{ME07}.
Moreover, for consistent variable selection with the Lasso,
it is known that the so-called
``neighborhood stability condition''  (\cite{MB06}) for the design matrix,
which has been re-formulated in a nicer form as the ``irrepresentable condition''
(\cite{ZY07}),
is sufficient and essentially necessary.
 \cite{Wai07,Wai08}
analyzes the smallest sample size needed
 to recover a sparse signal under certain incoherence conditions,
Because   irrepresentable or  incoherence conditions
are restrictive and much stronger than 
restricted eigenvalue conditions (see 
\cite{vdG:2009} for a comparison), we conclude that the Lasso 
for exact variable selection only works in a
rather narrow range of problems, excluding for example some cases where the design
exhibits strong (empirical) correlations.
 
Regularization with the $\ell_q$-``norm" with $q <1$ will mitigate
some of the bias problems, see \cite{zhang2010nearly}.
Related are
multi-step procedures where each of the steps involves a convex
optimization only. A prime example is the adaptive Lasso
which is a 
two-step algorithm and whose repeated application corresponds in some
``loose'' sense to a non-convex penalization scheme (\cite{zouli2008}). 
 \cite{Zou06} analyzed the adaptive Lasso in an asymptotic setup for the case where $p$ is
fixed. Further progress in the
high-dimensional scenario has been achieved 
by \cite{HMZ08}. 
Under a rather strong mutual incoherence condition between 
every pair of relevant and irrelevant covariables, they prove that the adaptive Lasso recovers the correct  
model and has an oracle property.
As we will explain in Subsection \ref{irrepresentable.section}, the adaptive Lasso 
indeed essentially needs a - still quite restrictive - weighted version of the 
irrepresentable condition in order to be able to correctly estimate the support
of the coefficients.

\cite{MY08} examine the thresholding procedure, assuming all non-zero components are large enough, an assumption we will avoid.
Thresholding and multistage procedures are also considered in
\cite{candes2006stable}, \cite{candes14enhancing}.
In \cite {Zhou09, Zhou10}, it is shown that a multi-step
thresholding procedure can accurately estimate a sparse vector $\beta \in \R^p$
under the restricted eigenvalue condition of \cite{bickel2009sal}.
The two-stage procedure in \cite{Zhang09} applies ``selective penalization" in the
second stage. This procedure is studied assuming incoherence conditions.
A more general framework for multi-stage variable selection 
was studied by \cite{WR08}. Their
approach controls the probability of false positives (type I error) but
pays a price in terms of false negatives (type II error).
The main contribution of this paper is that we provide bounds for the adaptive Lasso
that are comparable to the bounds for 
the Lasso followed
by a thresholding procedure. Because the true regression itself,
or its linear projection, 
is perhaps not sparse, we moreover consider a sparse
approximation of the truth, somewhat in the spirit of \cite{zhang2008sparsity}.

\subsection{Organization of the paper}
The next section introduces the sparse oracle approximation, with which we compare the initial and
adaptive Lasso.
In Section \ref{preview.section}, we present  the main results.
Eigenvalues and their
restricted and sparse counterparts are defined in Section \ref{notations.section}.
Some conclusions are presented in Section \ref{conclusion.section}.

The rest of the paper presents intermediate results and
complements for establishing the main results of Section \ref{preview.section}.
In Section \ref{noiseless.section}, we consider the noiseless
case, i.e., the case where $\epsilon = 0 $. The reason is that many of the
theoretical issues involved concern the approximation properties of the two stage procedure, and not so much the fact that there is noise.
 By studying the
noiseless case first, we separate the approximation problem from the
stochastic problem.

Both initial and adaptive Lasso are special cases of a weighted Lasso.
We discuss prediction error, $\ell_q$-error ($q\in \{ 1, 2 \} $) and variable selection with the
weighted Lasso 
in Subsection \ref{weight.section}. Theorem \ref{weight} in this section is the core of the present work, as regards prediction and estimation.
Lemma \ref{weightselect} in this section is the main result as regards variable selection.
 The behavior of the noiseless initial and adaptive Lasso are simple
corollaries of Theorem \ref{weight} and Lemma \ref{weightselect}.
We give in Subsection \ref{init.section}
the resulting bounds for the initial Lasso and discuss in Section \ref{thres.section} its
thresholded version. In Subsection \ref{adap.section} we derive results for the adaptive
Lasso by comparing it with a thresholded
initial Lasso. 
Moreover,
Subsection \ref{irrepresentable.section} briefly discusses the 
weighted irrepresentable condition, to show that even  the adaptive Lasso
needs strong conditions on the design
for exact variable selection. This subsection is linked to Corollary
\ref{min.corollary}, where it is proved that the false positives of
the adaptive Lasso vanish if the coefficients of the oracle are sufficiently large.

Section \ref{noise.section} studies the noisy case. 
It is an easy extension of the results of Sections \ref{weight.section},
\ref{init.section}, \ref{thres.section} and \ref{adap.section}. 
We do however need to further specify the choice of the tuning parameters
$\lambda_{\rm init}$ and $\lambda_{\rm adap}$.
After explaining the notation, we present the bounds for the
prediction error, estimation error and for the number
of false positives,  of the weighted Lasso.
This then provides us with the tools to prove the main results.

 All  proofs are in Section \ref{proofs.section}.
Here, we also present  explicit constants in the bounds to highlight
the non-asymptotic character of the results. 

\section{Model misspecification, weak variables and the oracle}\label{target}
Let
$$\EE {\bf Y} := {\bf f}^0 , $$ 
where ${\bf f}^0$ is the regression function. 
First, we note that without loss of generality, we can assume that ${\bf f}^0$ is
linear.
If ${\bf f}^0$ is non-linear in the
covariates, we consider its projection 
${\bf X} \beta_{\rm true}$
onto the linear space
$\{ {\bf X} \beta : \beta \in \R^p \}$,
i.e.,
$$ {\bf X} \beta_{\rm true} := \arg \min_{{\bf X} \beta }
\| {\bf f}^0 - {\bf X} \beta\|_2 . $$
It is not difficult to see that  all our results 
still hold if ${\bf f}^0$ is replaced by its projection ${\bf X} \beta_{\rm true}$. The statistical implication is very relevant.
The mathematical argument is the orthogonality
$${\bf X}^T ( {\bf X} \beta_{\rm true} - {\bf f}^0 ) = 0 . $$
  For ease of notation, we therefore assume from now on that ${\bf f}^0$ is indeed
linear: 
$${\bf f^0}: = {\bf X} \beta_{ {\rm true}} .$$

Nevertheless, $\beta_{\rm true}$ itself may not be sparse.
Denote the active set of $\beta_{\rm true}$ by
$$S_{\beta} := \{ j : \ \beta_{j, {\rm true}} \not= 0 \}, $$
which has cardinality
$s_{\rm true} := | S_{\rm true}|$. It may well be that $s_{\rm true}$ is
quite large, but that there are many weak variables, that is, many very small non-zero coefficients in
$\beta_{\rm true}$.
 Therefore, the sparse object we aim to recover may not be the ``true" unknown parameter $\beta_{\rm true} \in \R^p$ of the linear regression, but rather a sparse approximation.
We believe that an extension to the case where ${\bf f}^0$ is only
``approximately" sparse,  better reflects the true state of nature.
We emphasize however that throughout the paper, it is allowed to replace the
oracle approximation $b^0 $ given below by $\beta_{\rm true}$. This would
simplify the theory. However, we have chosen not to follow this
route because it generally leads to a large price to pay in the bounds.

The sparse approximation of ${\bf f}^0$ that we consider 
 is defined as follows.
For a set of indices $S \subset \{ 1 , \ldots , p \}$ and for $\beta \in \R^p$, we let
$$\beta_{j,S} := \beta_j {\rm l} \{ j \in S\} , \ j=1 , \ldots , p . $$
Given a set $S $,
the best approximation of ${\bf f}^0$ using only variables in $S$ is
$${\rm f}_S= {\bf X} b^S:= \arg \min_{f = {\bf X} {\beta_S}} \| f - {\bf f}^0 \|_2 . $$

Thus, ${\rm f}_S$
 is the projection of ${\bf f}^0$ on the linear span of the variables in $S$.
Our target is now 
the projection ${\rm f}_{S_0}$, where
$$S_0 := \arg \min_{S \subset S_{\rm true} }  \biggl \{
\| {\rm f}_S - {\bf f}^0 \|_2^2 /n + 7 \lambda_{\rm init}^2 | S| / \phi^2 (6,S) \biggr \} . $$
Here, $|S|$ denotes the size of $S$. Moreover, $\phi^2 (6,S)$ is
a ``restricted eigenvalue" (see Section \ref{notations.section} for its definition), 
which depends on the Gram matrix $\hat \Sigma$ and on the set $S$.
The constants
are chosen in relation with the oracle result (see Corollary \ref{noisyinit.corollary}). 
In other words, ${\rm f}_{S_0}$ is the optimal $\ell_0$-penalized
approximation, albeit that it is discounted by the restricted eigenvalue
$\phi^2 (6, S_0)$. 
To facilitate the interpretation, we require $S_0 $ to be a subset of
$S_{\rm true}$, so that the oracle is not allowed to trade irrelevant
coefficients against restricted eigenvalues.
With $S_0 \subset S_{\rm true}$,
any false positive selection with respect to  $S_{\rm true}$ is also a false positive for $S_0$.

We refer to ${\rm f}_{S_0}$ as the ``oracle".
The set $S_0$ is called the oracle active set, and 
$b^0 = b^{S_0}$ are the oracle coefficients, i.e., 
$${\rm f}_{S_0} = {\bf X} b^0 .$$
We write $s_0 = | S_0| $.

Inferring the sparsity pattern, i.e. variable selection, 
refers to the task of estimating the set of non-zero
coefficients,  that is,  to have a limited number
of false positives (type I errors) and false negatives 
(type II errors).
It can be verified that under reasonable conditions
with suitably chosen tuning parameter
$\lambda$, 
the ``ideal" estimator
$$\hat \beta_{\rm ideal} := \arg \min_{\beta}
\biggl \{ \| {\bf Y} - {\bf X} \beta \|_2^2 /n + \lambda^2 
| \{ j: \ \beta_j \not= 0 \} | \biggr \} , $$
has $O( \lambda^2 s_0)$ prediction error and
$O(s_0)$ false positives (see for instance \cite {Barr:Birg:Mass:1999}
and \cite{vdG:2001}).
With this in mind, we generally aim at $O(s_0)$ false positives
(see also
\cite{Zhou10}), yet keeping the
prediction error as small as possible
(see Corollary \ref{compare.corollary}).

As regards false negative selections, we refer to Subsection \ref{falsenegatives.section}, where
we derive bounds based on the $\ell_q$-error.

\section{Main results} \label{preview.section}

\subsection{Main conditions}

The behavior of the thresholded Lasso and adaptive Lasso depends on the tuning parameters, on the design, as well as
on the true ${\bf f}^0$, and actually on the interplay between these quantities.
To keep the exposition clear,
we will use order symbols. Our expressions are functions of $n$, $p$,
${\bf X}$, and ${\bf f}^0$, and also of the tuning parameters 
$\lambda_{\rm init}$, $\lambda_{\rm thres}$,  and $\lambda_{\rm adap}$.
For positive functions $g$ and $h$, we say that $g=O(h)$ if $
\| g / h\|_{\infty} $ is bounded, and
$g \asymp h$ if in addition $\| h/g\|_{\infty} $ is bounded.
Moreover, we say that $g = O_{\rm suff }(h)$ if
$\| g / h \|_{\infty} $ is not larger than a suitably chosen 
sufficiently small constant,
and
$g \asymp_{\rm suff} h $ if in addition $\| h / g \|_{\infty} $ is bounded. 

Our results depend on restricted eigenvalues
$\phi(L, S, N)$, minimal restricted eigenvalues $\phi_{\rm min} (L, S, N)$, and
minimal sparse eigenvalues $\phi_{\rm sparse} ( S , N)$
(which we generally think of as being not too small),
as well on maximal sparse eigenvalues $\Lambda_{\rm sparse} (s)$
(which we generally think of being not too large).
The exact definition of these constants is given in  
Section \ref{notations.section}.  

To simplify the expressions, we assume throughout that 
\begin{equation} \label{biasvariance}
\| {\rm f}_{S_0} - {\bf f}^0 \|_2^2 / n = O( \lambda_{\rm init}^2
s_0 / \phi^2 (6, S_0))
\end{equation}
 (where $\phi (6, S_0) = \phi( 6 , S_0 , s_0) $),
which roughly says that the oracle ``squared bias" term is not substantially
larger than the oracle ``variance" term. 
For example, in the case of orthogonal design, this condition holds if
the small non-zero coefficients are small enough, or if there are not too many of them,
i.e., if
$$\sum_{|\beta_{j, {\rm true}} |^2 \le 7 \lambda_{\rm init}^2 }
|\beta_{j,{\rm true} } |^2 = O (\lambda_{\rm init}^2 s_0 ) . $$
We stress that (\ref{biasvariance}) is merely to write order bounds for
the oracle, bounds with which we compare the ones for the various Lasso
versions. If actually the ``squared bias" term is the dominating term, this mathematically
does not alter the theory but makes the result more difficult to
interpret. 

We will furthermore discuss the results on the set
$${\cal T} := \biggl \{ 4 \max_{1 \le j \le p} | \epsilon^T {\bf X}_j / n | \le \lambda_{\rm init}
\biggr  \} ,$$
where ${\bf X}_j$ is the $j$-th column of the matrix ${\bf X}$.
For an appropriate choice of $\lambda_{\rm init}$, 
depending on the distribution of $\epsilon$, the set ${\cal T}$ has
large probability. Typically, $\lambda_{\rm init}$ can be taken of order
$$\sqrt { \log p / n } .$$
The next lemma serves as an example, but the
results can clearly be extended to other distributions.

\begin{lemma} \label{noise.lemma}
Suppose that $\epsilon \sim {\cal N} (0, \sigma^2 I)$. 
Take for a given $t>0$, 
$$\lambda_{\rm init} = 4 \sigma \sqrt { 2t + 2\log p \over n } .$$
Then
$$\PP ( {\cal T} ) \ge 1- 2 \exp [-t] . $$
\end{lemma}
%Extension to other distributions is straightforward, as we show in Section \ref{noise.section}. The results may furthermore be improved when the largest eigenvalue
%of the Gram matrix $\hat \Sigma = {\bf X}^T {\bf X}/n$ is well-behaved (e.g., when it
%is $O(1)$, see Section \ref{noise.section}).

The following conditions play an important role.
Conditions A and AA for thresholding are similar to those in \cite{Zhou10}
(Theorems 1.2, 1.3 and 1.4). 

{\bf Condition A} {\it For the thresholded Lasso, the threshold level $\lambda_{\rm thres}$ is chosen sufficiently large,  in such
a way that
$$\biggl [ {1 \over \phi^2(6, S_0 ,  2 s_0) }\biggr ] 
\lambda_{\rm init}= O_{\rm suff} (\lambda_{\rm thres}) . $$}

{\bf Condition AA} {\it For the thresholded Lasso, the threshold level $\lambda_{\rm thres}$ is chosen sufficiently large, but
such that
$$\biggl [ {1 \over \phi^2(6, S_0 ,  2 s_0) } \biggr ] \lambda_{\rm init}\asymp_{\rm suff}  \lambda_{\rm thres} . $$}

{\bf Condition B} {\it For the adaptive Lasso, the tuning parameter $\lambda_{\rm adap}$ is chosen
sufficiently large, in such a way that
$$ \biggl [ {\Lambda_{\rm sparse}
(s_0) 
 \over \phi_{\rm min}^3 ( 2 , S_0 , 2 s_0) }\bigg ] \lambda_{\rm init}= O_{\rm suff} (\lambda_{\rm adap}) . $$}
 
 {\bf Condition BB} {\it For the adaptive Lasso, the tuning parameter $\lambda_{\rm adap}$ is chosen
sufficiently large, but such that
$$\biggl [  { \Lambda_{\rm sparse}
(s_0) 
 \over \phi_{\rm min}^3 ( 6 , S_0 , 2 s_0) }\biggr ] \lambda_{\rm init}  \asymp_{\rm suff} \lambda_{\rm adap} . $$}
 
% {\bf Condition C} {\it The tuning parameter $\lambda_{\rm adap}$ is chosen at least as
% large as the threshold $\delta$.}

The above conditions can be considered with a zoomed-out look,
neglecting the expressions in the square brackets
($[ \cdots ]$), 
and a zoomed-in look, taking into account what is inside
the square brackets.
One may think of $\lambda_{\rm init}$ as the noise level (see e.g.\
Lemma \ref{noise.lemma},  with the $\log p $-term the price for not knowing the
relevant coefficients a priori).
Zooming out, Conditions A and B say that the threshold level $\lambda_{\rm thres}$ and
the tuning parameter $\lambda_{\rm adap}$ are required to be at least of the
same order as $\lambda_{\rm init}$, i.e., they should not drop
below the noise level. 
Assumption AA and BB put these parameters exactly
at the noise level, i.e., at the smallest value we allow. The reason to do this
is that one then can have good prediction and estimation bounds.
If we zoom in, we see in the square brackets the role played by the
various eigenvalues. As they are defined only later in Section \ref{notations.section},
it is at first reading perhaps easiest to remember that the $\phi$'s can be
small and the $\Lambda$'s can be large, but one hopes they behave well,
in the sense that the values in the square brackets are not too large.

\subsection{The results}

The next three theorems contain the  main ingredients of the present work.
Theorem \ref{noisyinit.theorem} is not new (see e.g.\ \cite{Bunea:06, Bunea:07a, Bunea:07b}, \cite{bickel2009sal},
\cite{koltch09a}), albeit that we replace the perhaps non-sparse
$\beta_{\rm true}$ by the sparser $b^0$ (see also
\cite{vandeG08}). Recall that the latter replacement is done because it yields generally an improvement
of the bounds.

\begin{theorem}\label{noisyinit.theorem}
For the initial Lasso $\hat \beta_{\rm init} = \hat \beta$ defined in (\ref{eq.init}),
we have on ${\cal T}$, 
$$\| {\bf X} \hat \beta_{\rm init}  - {\bf f}^0 \|_2^2 / n =\biggl [
{1 \over \phi^2 ( 6, S_0 ) }\biggr]  O( \lambda_{\rm init}^2 s_0 ) , $$
and
$$\| \hat \beta_{\rm init} - b^0 \|_1  =\biggl [  {1 \over \phi^2 ( 6, S_0) }\biggr ] O( \lambda_{\rm init} s_0 ) , $$
and
$$\| \hat \beta_{\rm init} - b^0 \|_2= 
\biggl [ {1 \over \phi^2 ( 6, S_0 , 2 s_0) }\biggr ] O( \lambda_{\rm init} \sqrt {s_0} ) .$$
\end{theorem}

The next theorem discusses thresholding. The results correspond to those in
\cite{Zhou10}, and will be invoked to prove similar bounds for the adaptive
Lasso, as presented in Theorem \ref{noisyadap.theorem}.

\begin{theorem} \label{noisythres.theorem}
Suppose Condition A holds. Then on ${\cal T}$,
$$\| {\bf X} \hat \beta_{\rm thres} - {\bf f}^0 \|_2^2 /n  = 
  \biggl [ 
 { \Lambda_{\rm sparse}^2 (s_0) } \biggr ] 
 {\lambda_{\rm thres}^2 \over \lambda_{\rm init}^2 } O (\lambda_{\rm init}^2 s_0)   , $$
 and
 $$\| \hat b_{\rm thres} - b^0 \|_1  =
  \biggl [
 {\Lambda_{\rm sparse} (s_0) \over\phi_{\rm sparse}
 ( S_0, 2s_0 ) }  \biggr ] {\lambda_{\rm thres} \over \lambda_{\rm init} }
 O( \lambda_{\rm init}  {s_0} ) , $$
 and
 $$\| \hat b_{\rm thres} - b^0 \|_2  =
  \biggl [
 {\Lambda_{\rm sparse} (s_0) \over\phi_{\rm sparse}
 ( S_0, 2s_0 ) }  \biggr ] { \lambda_{\rm thres}  \over  \lambda_{\rm init}} O( \lambda_{\rm init} \sqrt {s_0} ) , $$
 and
 $$| \hat S_{\rm thres} \backslash S_0 | = \biggl [
 { 1 \over \phi^4 ( 6 , S_0 , 2s_0) } \biggr ]  {\lambda_{\rm init}^2 \over \lambda_{\rm thres}^2} 
 O(s_0) . $$

 \end{theorem}

 \begin{theorem} \label{noisyadap.theorem}
 Suppose Condition B holds. Then on ${\cal T}$, 
$$
 \| {\bf X} \hat \beta_{\rm adap} - {\bf f}^0 \|_2^2 / n  = \biggl [ { \Lambda_{\rm sparse} (s_0)
  \over
 \phi_{\rm min} ( 6 , S_0 , 2 s_0) }  \biggr ] {\lambda_{\rm adap}
 \over \lambda_{\rm init} } O ( \lambda_{\rm init}^2 
 s_0 ) , 
$$
 and
 $$ \|\hat \beta_{\rm adap} - b^0 \|_1 =
 \biggl [ { \Lambda_{\rm sparse}^{1/2} (s_0) \over
 \phi_{\rm min}^{3/2} ( 6 , S_0 , 2 s_0) }  \biggr ] \sqrt {  \lambda_{\rm adap} \over
 \lambda_{ \rm init}}
 O ( \lambda_{\rm init} 
s_0 ) , $$
 and
  $$ \|\hat \beta_{\rm adap} - b^0 \|_2 =
 \biggl [ { \Lambda_{\rm sparse}^{1/2}  (s_0) 
 \phi_{\rm min}^{1/2} (6 , S_0 ,2 s_0) \over
 \phi_{\rm min}^2 ( 6 , S_0 , 3 s_0) }  \biggr ]\sqrt {  \lambda_{\rm adap} \over
 \lambda_{ \rm init}}
O( \lambda_{\rm init} 
 \sqrt {s_0 } ) , $$
 and
 $$| \hat S_{\rm adap} \backslash S_0|  = 
 \biggl [  {  \Lambda_{\rm sparse}^2 (s_0) \over
\phi^4 ( 6 , S_0 , 2 s_0) }  {\Lambda_{\rm sparse} (s_0) \over
\phi_{\rm min} ( 6 , S_0 , 2 s_0) } \biggr ] {\lambda_{\rm init} 
\over \lambda_{\rm adap} } O ( s_0) .  $$
\end{theorem}

We did not present a bound for the number of false positives of
the initial Lasso: it can be quite large depending on further conditions as given in
Lemma \ref{refineselect.lemma}. A rough bound is presented in Lemma
\ref{Lambdamax.lemma}. 

Theorem \ref{noisythres.theorem} and \ref{noisyadap.theorem}
show how the results depend on the choice of the tuning parameters
$\lambda_{\rm thres}$ and $\lambda_{\rm adap}$. The following corollary takes the
choices of Conditions AA and BB, as these choices give
the smallest prediction and estimation error.

\begin{corollary}\label{compare.corollary}
Suppose we are on ${\cal T}$. Then, 
under Condition AA,
 \begin{equation}\label{threspredict}
  \| {\bf X} \hat b_{\rm thres} - {\bf f}^0 \|_2^2 / n
  =\biggl [  { \Lambda_{\rm sparse}^2 (s_0)
 \over \phi^4 ( 6, S_0 , 2 s_0) } \biggr ] O ( \lambda_{\rm init}^2 s_0) , 
 \end{equation}
 and
  $$ \| \hat b_{\rm thres} - b^0 \|_1=
  \biggl [  { \Lambda_{\rm sparse} (s_0) \over \phi_{\rm sparse} ( S_0 , 2 s_0) \phi^2 ( 6 , S_0 , 2 s_0) } \biggr ]O ( \lambda_{\rm init}  {s_0} )  , $$
 and
 $$ \| \hat b_{\rm thres} - b^0 \|_2= 
 \biggl [ { \Lambda_{\rm sparse} (s_0) \over \phi_{\rm sparse} ( S_0 , 2 s_0) \phi^2 ( 6 , S_0 , 2 s_0) }\biggr ]  O ( \lambda_{\rm init} \sqrt {s_0} ) , $$
 and
 \begin{equation} \label{thresselect}
 | \hat S_{\rm thres} \backslash S_0 | = O (s_0) . 
  \end{equation}
  Similarly, under Condition BB, 
 \begin{equation} \label{adappredict}
  \| {\bf X} \hat \beta_{\rm adap} - {\bf f}^0 \|_2^2 / n
 = \biggl [ { \Lambda_{\rm sparse}^2 (s_0) \over
 \phi_{\rm min}^4 ( 6 , S_0 , 2 s_0) }  \biggr ] O ( \lambda_{\rm init}^2 
s_0 ) , 
 \end{equation}
 and
 $$ \|\hat \beta_{\rm adap} - b^0 \|_1 =
 \biggl [ { \Lambda_{\rm sparse} (s_0) \over
 \phi_{\rm min}^3 ( 6 , S_0 , 2 s_0) }  \biggr ] O (\lambda_{\rm init} 
 s_0 ) , $$
 and
 $$ \|\hat \beta_{\rm adap} -b^0 \|_2 =
 \biggl [ { \Lambda_{\rm sparse} (s_0)
 \over
 \phi_{\rm min}^2 ( 6 , S_0 , 3 s_0) \phi_{\rm min} (6 , S_0 , 2 s_0)  }  \biggr ] O ( \lambda_{\rm init} 
\sqrt {s_0 } ) , $$
 and 
 \begin{equation}\label{adapselect}
 | \hat S_{\rm adap} \backslash S_0|  = \biggl [ 
{  \Lambda_{\rm sparse}^2 (s_0) \phi_{\rm min}^2 ( 6, S_0 , 2 s_0) \over
\phi^4 ( 6 , S_0 , 2 s_0) } \biggr ]  O ( s_0) .  
\end{equation}

\end{corollary}

\begin{remark}
{\rm Note that our conditions on $\lambda_{\rm thres}$ and $\lambda_{\rm adap}$
depend on the $\phi$'s and $\Lambda$'s, which are unknown.
Indeed, our study is of theoretical nature, revealing common features
of thresholding and the adaptive Lasso.
Furthermore, it is possible to remove the dependence of the
$\phi$'s and $\Lambda$'s, when one imposes stronger sparse
eigenvalue conditions, along the lines of \cite{zhang2008sparsity}.
In practice, the tuning parameters
are generally chosen by cross validation. }
\end{remark}

\subsection{Comparison with the Lasso}

At the zoomed-out level, where all $\phi$'s and $\Lambda$'s are neglected,
we see that the thresholded Lasso (under Condition AA) and the
adaptive Lasso (under Condition BB) achieve the same order of magnitude for
the prediction error as the initial, one-stage
Lasso discussed in Theorem \ref{noisyinit.theorem}. The same is true for
their estimation errors. 
Zooming in on the $\phi$'s and the $\Lambda$'s, their error bounds
are generally larger than for the initial Lasso.

For comparison in terms of false positives, we need a corresponding bound
for the initial Lasso. In the paper of \cite{zhang2008sparsity}, one can find
results that ensure that also for the initial Lasso, modulo
$\phi$'s and $\Lambda$'s, the
number of false positives is of order $s_0$.
However, this result requires rather
involved conditions which also
improve the bounds for the adaptive and thresholded Lasso.
We briefly address this refinement in Subsection \ref{refine.section},
imposing a condition of similar nature as the one used in
\cite{zhang2008sparsity}. Also under these stronger conditions, the general message remains that thresholding and 
the adaptive Lasso can have similar prediction and estimation error as
the initial Lasso, and are often far better as regards variable selection

In this section, we confine ourselves to the following lemma.
Here,  $\Lambda_{\rm max}^2$ is the largest eigenvalue of $\hat
\Sigma$, which can generally be quite large.

\begin{lemma}\label{Lambdamax.lemma}
On ${\cal T}$,
$$|\hat S_{\rm init} \backslash S_0 | \le \biggl [ {\Lambda_{\rm max}^2 
\over \phi^2 ( 6 , S_0)} \biggr ] O(s_0) . $$
\end{lemma}

\subsection{Comparison between adaptive and thresholded Lasso}

When zooming-out, we see that the adaptive and thresholded Lasso  have
bounds of the same order of magnitude, for prediction, estimation and variable selection. 

At the zoomed-in level, the adaptive and thresholded
Lasso also have very similar bounds for the prediction error
(compare (\ref{threspredict}) with (\ref{adappredict})) in terms
of the $\phi$'s and $\Lambda$'s. A
similar conclusion holds for their estimation error. 
We remark that our choice of Conditions AA and BB for the tuning parameters 
is motivated by the fact that according to our theory,
these give the smallest prediction and estimation errors.
It then turns out that the ``optimal'' errors of the two methods match
at a quite detailed level.
However, if we zoom-in even further and
look at the definition of $\phi_{\rm sparse}$, $\phi$, and $\phi_{\rm min}$
in Section \ref{notations.section}, it will show up that the bounds for the adaptive Lasso  
prediction and estimation error are (slightly) larger.

Regarding variable selection, 
at zoomed-out level the results are also comparable
(see (\ref{thresselect}) and (\ref{adapselect})).
Zooming-in on the the $\phi$'s and $\Lambda$'s, the adaptive
Lasso may have more false positives than the thresholded version. 
 
A conclusion is that at the zoomed-in level, the adaptive Lasso has less
favorable bounds as the refitted thresholded Lasso. However, these are
still only bounds, which are based on focussing on a direct comparison
between the two methods, and we may have lost the finer properties of the
adaptive  Lasso. Indeed, the non-explicitness of the adaptive Lasso makes
its analysis a 
non-trivial task. The adaptive Lasso is a quite popular practical method,
and we certainly do not advocate that it should always be replaced by
thresholding and refitting.

\subsection{Bounds for the number of false negatives}\label{falsenegatives.section}

The $\ell_q$-error has immediate
consequences for the number of false negatives: if for some estimator $\hat \beta$,
some target $b^0$, and some constant $\delta_q^{\rm upper}$
one has
$$\| \hat \beta - b^0 \|_q \le \delta_q^{\rm upper}  $$
then the number of undetected yet large coefficients cannot be very large,
in the sense that
$$| \{ j: \ \hat \beta_j = 0 , | b_{j}^0|  > \delta \} |^{1/q} \le 
{\delta_q^{\rm upper}  \over \delta }. $$

Therefore, on ${\cal T}$, for example
$$\biggl |\biggl  \{ j: \ \hat \beta_{j, {\rm init}}  = 0 , 
\biggl [ {1 \over \phi^2 ( 6 , S_0 , 2 s_0) } \biggr ] \sqrt {s_0} \lambda_{\rm init}  =
O_{\rm suff} (| b_{j}^0| )\biggr \} \biggr |
= 0 . $$
Similar bounds hold for the thresholded and the adaptive Lasso
(considering now, in terms of the $\phi$'s and $\Lambda$'s, somewhat larger
$|b_j^0|$).

One may argue that one should not aim at detecting
variables that the oracle considers as irrelevant.
Nevertheless, given an estimator $\hat \beta$, it is straightforward to bound $\| \hat \beta - \beta_{\rm true}\|_q$
in terms of $\| \hat \beta - b^0 \|_q$:  
apply the triangle inequality
$$\| \hat \beta - \beta_{\rm true} \|_q \le \| \hat \beta - b^0 \|_q + 
\| b^0 - \beta_{\rm true}\|_q . $$ 
Moreover, for $q=2$, one has the inequality
$$ \| b^0 - \beta_{\rm true} \|_2^2  \le { \| {\rm f}_{S_0} - {\bf f}^0 \|_2^2 \over
n  \Lambda_{\rm min}^2 (S_{\rm true}) } , $$
where $ \Lambda_{\rm min} (S)$ is the smallest eigenvalue
of the Gram matrix corresponding to the variables in $S$. 
One may verify that
$\phi(6, S_{\rm true}) \le \Lambda_{\rm min} (S_{\rm true})$. 
In other words, by choosing $\beta_{\rm true}$ as target
instead of $b^0$, does in our approach not lead to an improvement in
the bounds for
$\| \hat \beta - \beta_{\rm true} \|_2$.

\subsection{Having large coefficients}

Let us have a closer look at what conditions on the
size of the coefficients can bring us. We only discuss the
adaptive Lasso (thresholding again giving similar results, see also \cite{Zhou10}).

We define
$$| b^0 |_{\rm min} := \min_{j \in S_0 } | b_j^0 | . $$
Moreover, we let
$$|b^0|_{\rm harm}^{2}  := 
\biggl ( {1 \over s_0} \sum_{j \in S_0} {1 \over | b_j^0 |^2 } \biggr )^{-1} $$
be the harmonic mean of the squared coefficients.

{\bf Condition C} {\it For the adaptive Lasso, take $\lambda_{\rm adap}$ sufficiently large,
such that
$$|b^0|_{\rm harm} = O_{\rm suff} (\lambda_{\rm adap} ) . $$}

{\bf Condition CC} {\it For the adaptive Lasso, take $\lambda_{\rm adap}$ sufficiently large,
but such that
$$|b^0|_{\rm harm} \asymp_{\rm suff} \lambda_{\rm adap} . $$}

\begin{lemma}\label{min.lemma}
Suppose that  for some constant $\delta_{\infty}^{\rm upper} $,
on ${\cal T}$, 
$$\| \hat \beta_{\rm init} - b^0 \|_{\infty}  \le \delta_{\infty}^{\rm upper}. $$
Assume in addition that
$$| b^0 |_{\rm min}> 2  \delta_{\infty}^{\rm upper} . $$
Then under Condition C, 
$$\| {\bf X} \hat \beta_{\rm adap}^2 - {\bf f}^0 \|_2^2 / n  =\biggl [ { 1 \over 
 \phi^2 (6, S_0 )}\biggr ] { \lambda_{\rm adap}^2
\over | b^0 |_{\rm harm}^2 }
O( \lambda_{\rm init}^2  s_0 
 )  ,$$
and 
$$\| \hat \beta_{\rm adap} -b^0\|_1  =  \biggl [ { 1 \over 
 \phi^2 (6, S_0 )}\biggr ]
{ \lambda_{\rm adap}
\over | b^0 |_{\rm harm} }
O  (
   { \lambda_{\rm init}  s_0  }  )  ,$$
 and
$$\| \hat \beta_{\rm adap} -b^0 \|_2 
=  \biggl [ { 1 \over 
 \phi^2 (6, S_0 , 2 s_0)}\biggr ]
{ \lambda_{\rm adap}
\over | b^0 |_{\rm harm} }
O  (
   { \lambda_{\rm init}  \sqrt {s_0 }  }  )  ,$$
 and
 $$| \hat S_{\rm adap}   \backslash S_0 | =\biggl [ 
{  \Lambda_{\rm sparse}^2 (s_0) \over  \phi^2 (6, S_0 )\phi^4 ( 6 , S_0 , 2 s_0) } \biggr ] 
 O \biggl (   { \lambda_{\rm init}^2   s_0  \over
  |b^0|_{\rm harm}^2 }\biggr ) .
$$
\end{lemma}

 It is clear that by Theorem \ref{noisyinit.theorem}, 
$$\| \hat \beta_{\rm init} - b^0 \|_{\infty} = 
  \biggl [  { \sqrt {s_0}  \over  \phi^2(2, S_0)}  \wedge {1  \over\phi^2 (2, S_0 , 2 s_0)} \biggr ] 
O( 
 \lambda_{\rm init} \sqrt {s_0}  ) . $$
This can be improved under coherence conditions on the Gram matrix.
To simplify the exposition, we will not discuss such improvements in detail
(see \cite{Lou08}).

Under Condition CC, the bound for the prediction error and estimation error is again the smallest.
We moreover have the following corollary for the number of false positives.

\begin{corollary}\label{min.corollary} Assume the conditions of Lemma
\ref{min.lemma} and
$$  \phi^2 ( 6 , S_0 , 2 s_0)
 \lambda_{\rm init} \sqrt {s_0} = O(| b^0 |_{\rm harm}  ).$$
Then
on ${\cal T}$, 
$$| \hat S_{\rm adap}   \backslash S_0 | =\biggl [ 
{  \Lambda_{\rm sparse}^2 (s_0) \over  \phi^2 (6, S_0 ) } \biggr ] 
 O  ( 1) .
$$
\end{corollary}
 
 By assuming that $| b^0|_{\rm harm}$ is sufficiently large, that is, 
 $$\biggr [   {  \Lambda_{\rm sparse} (s_0) \over  \phi (6, S_0 )\phi^2 ( 6 , S_0 , 2 s_0) } \biggr ] 
  \lambda_{\rm init} \sqrt {s_0} = O_{\rm suff}(| b^0 |_{\rm harm}  ),$$
  one can
 bring $|\hat S_{\rm adap} \backslash S_0 | $ down to zero, i.e.,  no false
 positives. One may verify that this boils down to a situation where the
 weighted irrepresentable condition holds: see 
 Example \ref{separate.example} in Subsection \ref{irrepresentable.section}. 
 
 As discussed in Section \ref{falsenegatives.section}, large non-zero
 coefficients also lead to a small number or eventually zero false negative selections. 
 Therefore, the adaptive and thresholded Lasso are recovering the support
 of $S_0$ if all of its non-zero coefficients are sufficiently large
 (in absolute value), assuming much weaker conditions on the design
 than the (unweighted) irrepresentable condition, which is necessary for the Lasso.

\section{Notation and definition of generalized eigenvalues}\label{notations.section}
We reformulate the problem in 
$L_2 (Q)$, where
$Q$ is a generic probability measure on some space ${\cal X}$.
(This is somewhat more natural in the noiseless case, which we will
consider in Section \ref{noiseless.section}.)
Let
$\{ \psi_j \}_{j=1}^p \subset L_2 (Q)$ be a given dictionary. 
For $j=1 , \ldots , p$, the function $\psi_j$ will play the role of the $j$-th co-variable. 
The Gram matrix is
$$\Sigma := \int \psi^T \psi d Q , \ \psi:= ( \psi_1 , \ldots , \psi_p ) . $$
We assume that $\Sigma$ is normalized, i.e., that
$\int \psi_j^2 d Q =1$ for all $j$. In our final results, we will actually take
$\Sigma= \hat \Sigma$, the (empirical) Gram matrix corresponding to
fixed design.

Write
a linear function of the $\psi_j$ with coefficients $\beta\in \R^p $ as
$$f_{\beta} := \sum_{j=1}^p \psi_j \beta_j . $$
The  $L_2 (Q)$-norm is denoted by $\| \cdot \|$, so that
$$\| f_{\beta} \|^2 = \beta^T \Sigma \beta . $$

Recall that for an arbitrary $\beta \in \R^p$, and an
 arbitrary index set $S$, we use the notation
 $$\beta_{j,S}= \beta_j {\rm l} \{ j \in S \} . $$

We now present our notation for eigenvalues.
We also introduce restricted eigenvalues and sparse eigenvalues.

\subsection{Eigenvalues} \label{eigenvalues.section}
 
 The largest eigenvalue of $\Sigma$ is denoted by $\Lambda_{\rm max}^2$,
 i.e.,
 $$\Lambda_{\rm max}^2 := \max_{\| \beta \|_2=1} \beta^T \Sigma \beta . $$
 We will also need the largest eigenvalue of a submatrix  containing the
 inner products of variables in $S$:
 $$\Lambda_{\rm max}^2 (S) := 
 \max_{\| \beta_S \|_2=1} \beta_S^T \Sigma \beta_S . $$
 Its minimal eigenvalue is
 $$\Lambda_{\rm min}^2 (S) := 
 \min_{\| \beta_S \|_2=1} \beta_S^T \Sigma \beta_S . $$

  \subsection{Restricted eigenvalues}
A restricted eigenvalue is of similar nature as
 the minimal eigenvalue of $\Sigma$, but with the coefficients $\beta$ restricted to
certain subsets of $\R^p$. 
The restricted eigenvalue condition we impose corresponds
to the so-called  {\it adaptive} version as introduced
in \cite{vdG:2009}. It differs from the restricted eigenvalue
condition in \cite {bickel2009sal}
or \cite{koltch09a, koltch09b}. This is due to the fact that we
want to mimic the oracle ${\rm f}_{S_0}$, that is, do not choose
${\bf f}^0$ as target, so that we have to deal with
a bias term $\| {\rm f}_{S_0} - {\bf f}^0 \|$. 
For a given
$S$, our restricted eigenvalue condition is stronger than the one
in \cite {bickel2009sal}
or \cite{koltch09a, koltch09b}. On the other hand, we apply
it to the smaller set $S_0$ instead of to $S_{\rm true}$.

Define for an index set $S\subset \{ 1 , \ldots , p \}$, and for a set ${\cal N} \supset S$ and constant $L > 0$,  the sets of restrictions
$${\cal R} (L,S, {\cal N} ) := \biggl \{ \beta: \ \| \beta_{{\cal N}^c} \|_1 \le L \sqrt {|{\cal N}|} \| \beta_{\cal N}  \|_2 , \ \max_{j \in {\cal N}^c } | \beta_j | \le \min_{j \notin
{\cal N}\backslash S } | \beta_j | \biggr \} . $$

{\bf Definition: Restricted eigenvalue.}
{\it For $N \ge |S|$, we call
$$\phi^2 (L, S,N) :=  
\min \biggl \{ {
\| f_{\beta} \|^2  \over \| \beta_{\cal N} \|_2^2 } :\  {\cal N} \supset S ,\ | {\cal N} | \le N , \ 
\beta \in {\cal R} (L,S, {\cal N} )   \biggr \} 
 $$
the {\rm $(L,S,N)$-restricted eigenvalue}.
The {\rm $(L,S,N)$-restricted eigenvalue condition} holds if $\phi(L,S,N)
>0$. \\
For the case $N=|S|$, we write $\phi(L,S):= \phi(L,S,|S|)$.\\
The {\rm minimal $(L,S,N)$-restricted eigenvalue} is
$$\phi_{\rm min}^2 (L, S,N) :=  \min_{{\cal N} \supset S, \ | {\cal N} | = N }
\phi^2 (L, {\cal N} ) . $$
}

It is easy to see that $\phi_{\rm min} (L, S, N) \le
\phi(L, S , N) \le \phi (L,S) \le \Lambda_{\rm min} (S)$ for all $L > 0 $.
It can moreover be shown that
$$\phi^2 (L, S,2|S|) \ge \min \biggl \{ 
\| f_{\beta} \|^2:\  {\cal N} \supset S , \ |{\cal N} | = 2|S|   , \ \| \beta_{{\cal N}^c}  \|_2 \le 1 , \ \| \beta_{\cal N} \|_2 = 1  \biggr \} . $$

\subsection{Sparse eigenvalues}

The fact that we also need 
sparse eigenvalues is in line with  the sparse Riesz condition occurring in \cite{zhang2008sparsity}.

{\bf Definition: Sparse eigenvalues.}
{\it For $N \in \{ 1 , \ldots , p \}$, the {\rm maximal sparse eigenvalue} is
$$\Lambda_{\rm sparse} (N) = \max_{{\cal N} : \ |{\cal N}| =N } 
\Lambda_{\rm max} ({\cal N}) . 
$$
For an index set $S \subset \{ 1 , \ldots , p \}$ with $|S| \le N$, 
the {\rm minimal sparse eigenvalue} is
$$\phi_{\rm sparse}(S, N)  := \min_{{\cal N} \supset S:\  |{\cal N}| = N } \Lambda_{\rm min} ({\cal N}) . $$}

One easily verifies that for any set ${\cal N}$ with $| N| =ks$, $k \in \Nat$, 
$$\Lambda_{\rm max} ({\cal N} ) \le \sqrt {k} \Lambda_{\rm sparse} (s) . $$
Moreover, for all $L\ge 0$, 
$$\phi_{\rm sparse}  ( S , N) = \phi ( 0 , S, N) \ge \phi ( L , S , N) . $$

\section{Conclusions}\label{conclusion.section} 

We present some comparable bounds for the adaptive Lasso and the thresholded
Lasso with refitting and we also compared them to the ordinary Lasso. The
framework of our analysis allows for misspecified 
linear models whose best linear projection is not necessarily sparse and
with possibly small non-zero regression coefficients, i.e., many  weak
variables. This setting is much more realistic than the usual
high-dimensional framework where the model is true with only a few but
strong variables.  

 Estimating the support $S_0$ of the non-zero coefficients
 is a hard statistical problem. The irrepresentable condition, which is essentially
 a necessary condition for exact recovery of the non-zero coefficients by the
 one-step Lasso, is much too restrictive in many cases. In this paper, our main
 focus is on having $O(s_0)$ false positives while achieving good prediction
 and estimation. This is inspired by the behavior
 of the ``ideal" $\ell_0$-penalized estimator.  
 
 We have examined thresholding the Lasso with
 least squares refitting and the adaptive Lasso.
 Our main conclusion is that both methods can have about the
 same prediction and estimation error as the one-stage ordinary Lasso,
 and that both gain over the one-stage Lasso in the sense of
 having less false positives.
 Moreover, according to our theory
(and not exploiting the fact that the adaptive Lasso mimics
thresholding and refitting using an ``oracle" threshold), 
thresholding with least squares refitting and the
adaptive Lasso perform equally well, even when considered at a
rather fine scale. Our bounds for the adaptive Lasso
are  more sensitive to
small (minimal) restricted eigenvalues or small minimal sparse eigenvalues,
or large sparse maximal eigenvalues.
Both thresholded and adaptive Lasso
benefit from a situation with large non-zero coefficients of the oracle. 

We do not give an account of the tightness of our bounds. The
thresholded Lasso allows a rather direct analysis, and we believe there
is little room for improvement of the bounds for this method. The analysis
of the adaptive Lasso more involved. Our comparison to thresholding might not do justice
to the adaptive Lasso. Indeed, we have not fully exploited the finer
oracle properties of the adaptive Lasso.
 
In practice the
the tuning parameters are often chosen by cross validation, which may
correspond to a choice giving the smallest
 prediction error. It is not within the scope of this paper
 to prove that with cross validation, thresholding and
 the adaptive Lasso again have comparable theoretical
 performance, although we do believe this to be the case.
 As for the computational aspect, we observe the following.
For the solution path for all $\lambda_{\rm adap}$, the adaptive Lasso
needs $O(n|\hat S_{\rm init}| \min(n, |\hat S_{\rm init}|))$ essential operation counts.
The same order of operation counts is needed when computing the thresholded Lasso
for the whole solution path over all $\lambda_{\rm thres}$.
 Therefore, the two methods are also computationally
 comparable.

\section{The noiseless case}\label{noiseless.section}

 Consider a fixed target ${\bf f}^0 = f_{\beta_{\rm true}} \in L_2 (Q)$. 
  Let $S \subset \{ 1 , \ldots , p \}$ and let 
  ${\rm f}_S := \arg \min_{f = f_{\beta_S}}  \| f_{\beta_S } - {\bf f}^0 \|  $
be the projection of ${\bf f}^0$ on the
$|S|$-dimensional linear space spanned by the variables $\{ \psi_j \}_{j \in S} $. 
We denote the coefficients of ${\rm f}_S$ by $b^S$, i.e., 
$${\rm f}_S = \sum_{j \in S} \psi_j b_{j}^S = f_{b^S}. $$
 The oracle set $S_0$ is defined by trading off
dimension against fit, namely 
\begin{equation} \label{defineS_0}
S_0 := \arg \min_{S \subset S_{\rm true}}  \biggl \{ \| {\rm f}_S - {\bf f}^0 \|^2 +{ 3 \lambda_{\rm init}^2
|S|  \over \phi^2 (2,S) } \biggr \} , 
\end{equation}
where the constants are now from Theorem \ref{weight} (or its
Corollary \ref{init.corollary}). 
 We call ${\rm f}_{S_0}$ the oracle, and we let $b^0:= b^{S_0} $,
 i.e., ${\rm f}_{S_0} = f_{b^0} $.

 For simplicity, we assume throughout that
 $$\| {\rm f}_{S_0} - {\bf f}^0 \|^2 = O ( \lambda_{\rm init}^2 s_0 / \phi^2 (2, S_0) ),
$$
which roughly says that the approximation error does not overrule
the penalty term.

The initial Lasso is
$$
\beta_{\rm init} := \arg \min_{\beta} \biggl \{ \| f_{\beta} - {\bf f}^0 \|^2 +
 \lambda_{\rm init} \| \beta \|_1 \biggr \} . 
$$
We assume that the tuning parameter $\lambda_{\rm init}$
is set at some fixed value. Of course, in the noiseless case, the optimal -
in terms of prediction error - value for $\lambda_{\rm init}$ is
$\lambda_{\rm init} =0$. 
However, in the noisy case, a strictly positive lower bound for $\lambda_{\rm init}$
is dictated by the noise level.
Write
\begin{equation} \label{initdef}
f_{\rm init} := f_{\beta_{\rm init}} , \ S_{\rm init}:= \{ j:\ \beta_{j , {\rm init}} \not= 0 \} , 
 \ \delta_{\rm init} := \| f_{\rm init}- {\bf f}^0 \| . 
  \end{equation}

Let for $\delta > 0$,
$$S_{\rm init}^{\delta} := \{ j : \ | \beta_{j, {\rm init}} | > \delta \} . $$
Then ${\rm f}_{S_{\rm init}^{\delta}} =f_{b^{S_{\rm init}^{\delta}}} $ is
the refitted Lasso after thresholding at $\delta$.
Note that we express explicitly the dependence
of the thresholded estimator on the threshold level,
which we now call $\delta$ (instead of $\lambda_{\rm thres}$
as we did in the introduction). 
The reason for this is that the analysis of the adaptive Lasso
will go via the thresholded Lasso with a choice of the threshold
$\delta$ that trades off prediction error against estimation error
(see (\ref{choicedelta}) in the proof of Theorem \ref{adap.theorem}).

The adaptive Lasso is
$$\beta_{\rm adap} := \arg \min_{\beta} \left \{ \| f_{\beta} - {\bf f}^0 \|^2 +
 \lambda_{\rm init} \lambda_{\rm adap} \sum_{j=1}^p { |\beta_j | \over
| \beta_{j, {\rm init}} | }  \right \} . $$
The second stage tuning parameter $\lambda_{\rm adap}$ is again assumed to be strictly positive.
We denote the resulting adaptive variants of (\ref{initdef}) by
$$f_{\rm adap} := f_{\beta_{\rm adap}} , \ S_{\rm adap}:= \{ j:\ \beta_{j , {\rm adap}} \not= 0 \} , \ \delta_{\rm adap} :=
\| f_{\rm adap} - {\bf f}^0 \| . $$

As the initial and adaptive Lasso are special cases of the weighted Lasso,
many of the results in Subsections \ref{init.section}, \ref{thres.section} and \ref{adap.section} are consequences of those
for the weighted Lasso as studied in Subsection \ref{weight.section}.
The weighted Lasso is
$$\beta_{\rm weight} := \arg \min_{\beta} 
\left \{ \| f_{\beta} - {\bf f}^0 \|^2 + \lambda_{\rm init} \lambda_{\rm weight}
\sum_{j=1}^p w_j |\beta_j | 
\right \}  , $$
where the $\{ w_j \}_{j=1}^p $ are non-negative weights.

We set $f_{\rm weight} := f_{\beta_{\rm weight} } $,
$S_{\rm weight} := \{ j: \ \beta_{j, {\rm weight}} \not= 0 \}$.
Moreover, we define
$$\| w_S \|_2^2  := \sum_{j \in S} w_j^2 , \ 
w_{S^c}^{\rm min} := \min_{j \notin S} w_j . $$

By the reparametrization $\beta \mapsto \gamma := W \beta$, 
where $W= {\rm diag} ( w_1 , \cdots , w_p )$, one
sees that the weighted Lasso is a standard Lasso with Gram matrix
$$\Sigma_{\rm weight} := W^{-1} \Sigma W^{-1} . $$
We emphasize however that $\Sigma_{\rm weight}$ is generally not normalized,
i.e., generally ${\rm diag} (\Sigma_{\rm weight}) \not= I $.

\subsection{The weighted Lasso}\label{weight.section}
We first present a bound for the prediction and estimation error and then consider variable
selection.

\begin{theorem} \label{weight}
Let $S$ be an index set with
cardinality $s:= |S|$, satisfying for some constants $M\ge 0$ and $L> 0$,
$$   w_{S^c}^{\rm min}  \ge M/L,  \  \| w_S \|_2 / \sqrt {s} \le M . $$
Then for all $\beta$,
we have
$$\| f_{\rm weight} - {\bf f}^0 \|^2    \le
2 \| f_{\beta_S} - {\bf f}^0 \|^2 + { 6 \lambda_{\rm init}^2
\lambda_{\rm weight}^2  M^2 s \over \phi^2 (2L,S) } 
 .$$
Moreover, for all $\beta$, we have
$$  \sqrt {s} 
\| (\beta_{\rm weight})_S - \beta_S \|_2 + 
 \| ( \beta_{ {\rm weight}})_{S^c} \|_1/L 
 \le 
{3 \| f_{\beta_S} - {\bf f}^0 \|^2
\over \lambda_{\rm init} \lambda_{\rm weight} M}   +{ 3 \lambda_{\rm init}
\lambda_{\rm weight} M s \over \phi^2 (2L,S) }  .$$
Finally, it holds for all $\beta$, that
$$\| \beta_{\rm weight} - \beta_S \|_2 \le
 { 6L \| f_{\beta_S}  - {\bf f}^0 \|^2  \over \lambda_{\rm init}
\lambda_{\rm weight} M  \sqrt {s_0} }
 + {6L \lambda_{\rm init} \lambda_{\rm weight} M ( s + s_0) \over \phi^2 (2L, S, s+s_0) 
  \sqrt s_0 }
  . $$
\end{theorem}

We will apply the above theorem with $S$ the set of the smaller weights.

\begin{corollary} \label{ordered.corollary}
Fix some arbitrary $\delta >0$, and let 
$$S_{\rm weight}^{\delta} \supset
\{ j: \ w_j < 1/\delta \} , \ 
(S_{\rm  weight}^{\delta})^c \supset \{ j : \ w_j > 1/ \delta \} .$$
The indices $j$ with $w_j = 1/ \delta$ can be put in either $S_{\rm weight}^{\delta}$
or in its complement.
Suppose that for some $\alpha \ge 0$, 
$$|S_{\rm weight}^{\delta}\backslash  S_0 |  \le \alpha s_0 . $$
Taking $S= S_{\rm weight}^{\delta}$, $L=1$ and $M = 1/ \delta $ in  Theorem \ref{weight},
we get that for all $\beta$, 
$$\| f_{\rm weight} - {\bf f}^0 \|^2  
 \le
%2 \| f_{\beta_{S_{\rm weight}^{\delta}}} - {\bf f}^0 \|^2 + { 6 \lambda_{\rm init}^2
%\lambda_{\rm weight}^2 | S_{\rm weight}^{\delta} | \over \phi^2 (2,S_{\rm weight}^{\delta}) } \|
%w_{S_{\rm weight}^{\delta}} \|_2^2 $$
%$$ \le 
2 \| f_{\beta_{S_{\rm weight}^{\delta}}} - {\bf f}^0 \|^2 + { 6 \lambda_{\rm init}^2
\lambda_{\rm weight}^2 (1+\alpha) s_0  \over \delta^2 \phi_{\rm min}^2 (2,S_0, (1+\alpha) s_0) } .
 $$
 Moreover,
 $$\| \beta_{\rm weight} - \beta_{S_{\rm weight}^{\delta} }\|_1 \le
 {3 \delta \| f_{\beta_{S_{\rm weight}^{\delta}}} - {\bf f}^0 \|^2 \over
 \lambda_{\rm init} \lambda_{\rm weight}} +
 {3 \lambda_{\rm init} \lambda_{\rm weight} (1+\alpha) s_0 \over
 \delta \phi^2 ( 2 , S_0 , (1+ \alpha ) s_0 ) } , $$
 and
$$\| \beta_{\rm weight} - \beta_{S_{\rm weight}^{\delta}} \|_2 
\le { 6 \delta  \| f_{\beta_{S_{\rm weight}^{\delta}} }  - {\bf f}^0 \|^2 
\over \sqrt {s_0} \lambda_{\rm init} \lambda_{\rm weight} }  + { 6 \lambda_{\rm init} \lambda_{\rm weight}(2+ \alpha) \sqrt {s_0}  \over \delta \phi_{\rm min}^2 (2, S_0 ,
(2+ \alpha) s_0) }  . 
 $$
 In the case $\alpha=0$, one may replace in the last bound, $\phi_{\rm min}^2 ( 2 , S_0 , 
 (2+ \alpha ) s_0) = \phi_{\min} ( 2 , S_0 , 2 s_0) $ by
 $\phi( 2 , S_0 , 2 s_0) $.
\end{corollary}

Our next theme is variable selection.
The {\it Karush-Kuhn-Tucker} ({\it KKT}) conditions (see \cite{bertsimas1997introduction})
can be invoked to derive Lemma \ref{weightselect} below,
where we use the notation
$$\| (1 / w)_S\|_2^2 := \sum_{j \in S} {1 \over w_j^2 } . $$

\begin{lemma}\label{weightselect} It holds that
\begin{equation} \label{stillin}
| S_{\rm weight} \backslash S_0 |^2 \le 4
\Lambda_{\rm max}^2 ( S_{\rm weight} \backslash S_0 ) 
{ \| f_{\rm weight} - {\bf f}^0 \|^2
\over \lambda_{\rm weight}^2 } {  \| (1 / w)_{S_{\rm weight}\backslash S_0 } \|_2^2
\over \lambda_{\rm init}^2 }  .
\end{equation}
If $| S_{\rm weight} \backslash S_0 |> s_0$, we have
$$| S_{\rm weight} \backslash S_0 | \le 
8\Lambda_{\rm sparse }^2 (s_0) {\| f_{\rm weight} - {\bf f}^0 \|^2   \over \lambda_{\rm weight}^2 s_0} 
{  \|( 1/w)_{S_{\rm weight} \backslash S_0}\|_2^2 \over \lambda_{\rm init}^2}.
$$ 
\end{lemma}

\subsection{The initial Lasso}\label{init.section}

Recall that
$$\delta_{\rm init}:= \| f_{\rm init} - {\bf f}^0 \| . $$
For $q \ge 1 $, we define
$$  \delta_q :=\| \beta_{\rm init} - b^0 \|_q .$$

\begin{theorem}\label{init.theorem} The prediction error of the
initial Lasso has
$$\delta_{\rm init}^2 = \biggl[ {1 \over  \phi^2 (2, S_0)} \biggr ]
O( \lambda_{\rm init}^2 s_0  ), $$
and its estimation error has
$$\delta_1 =  \biggl [ {1 \over  \phi^2 (2, S_0)} \biggr ] O ( \lambda_{\rm init} s_0 ), \ 
\delta_2 =  \biggl[  {1 \over  \phi^2 (2, S_0, 2 s_0 )} \biggr ]
O( \lambda_{\rm init} \sqrt {s_0} )  . $$
The initial estimator has number of false positives
$$| S_{\rm init} \backslash S_0 | =\biggl [  { \Lambda_{\rm max}^2
(S_{\rm init} \backslash S_0 )  \over \phi^2 (2, S_0) }
\biggr ] 
   O( s_0 ). $$
\end{theorem}

Considering the variable selection result, it is clear
that $\Lambda_{\rm max}^2 ( S_{\rm init} \backslash S_0 ) \le 
\Lambda_{\rm max}^2$. 
Without further conditions, this cannot be
refined, and the eigenvalue $\Lambda_{\rm max}^2$ can be quite large (yet having the minimal eigenvalue of $\Sigma$ bounded away from zero). Therefore, the result of  Theorem 
  \ref{init.theorem}
 needs further conditions for good variable selection properties of the initial Lasso.

\subsection{Thresholding the initial estimator}\label{thres.section}

Variable selection results by thresholding are not difficult to obtain:
$$
| S_{\rm init}^{\delta} \backslash S_0 |^{1/q} \le { \delta_q \over  \delta}  .
$$
Hence, for 
$\delta \ge 
\delta_1 / s_0 \wedge \delta_2 /\sqrt {s_0}$, we get for $q \in \{ 1 , 2 \}$, 
\begin{equation}\label{s0}
| S_{\rm init}^{\delta} \backslash S_0 |  \le s_0.
\end{equation}

If the coefficients of the oracle are sufficiently large, thresholding will
improve the prediction and estimation error.
Here, we do not impose such minimal size conditions.
The estimation error of the thresholded Lasso is then
still easy to assess. Our bound for the prediction error, however,
now depends on maximal sparse eigenvalues.

At this stage, we invoke the noiseless counterparts of Conditions A and AA.

{\bf Condition a} {\it We have
 $\lambda_{\rm init} / \phi^2 ( 2 , S_0) = O_{\rm suff} ( \delta) $.}
 
  {\bf Condition aa} {\it We have 
 $  \lambda_{\rm init} / 
 \phi^2 ( 2 , S_0 , 2 s_0)  \asymp_{\rm suff} \delta $.}

\begin{theorem}\label{thres.theorem} Assume Condition a. Then
$$\| {\rm f}_{S_{\rm init}^{\delta} } - {\bf f}^0 \|^2 =  \Lambda_{\rm sparse}^2 (s_0)
  \biggl [ 
 {\delta^2 \over \lambda_{\rm init}^2 } \biggr ] O (\lambda_{\rm init}^2 s_0)   , $$
 $$\| b^{S_{\rm init}^{\delta}} - b^0 \|_2  =
  {\Lambda_{\rm sparse} (s_0) 
 \over \phi_{\rm sparse}
 ( S_0 , 2s_0 ) } 
  \biggl [ 
 {\delta \over \lambda_{\rm init} }  \biggr ] O( \lambda_{\rm init} \sqrt {s_0} ) , $$
 and
 $$| S_{\rm init}^{\delta} \backslash S_0 | = \biggl [
 { 1 \over  \phi^4 ( 2 , S_0 , 2s_0) } \biggr ] 
 \biggl [
 { \lambda_{\rm init}^2 \over \delta^2 } \biggr ] 
 O(s_0) . $$
 The expressions for the prediction and estimation error lead to favoring the choice $ \lambda_{\rm init} / 
 \phi^2 ( 2 , S_0 , 2 s_0)  \asymp_{\rm suff} \delta $
  of Condition aa, which yields
 $$\| {\rm f}_{S_{\rm init}^{\delta} } - {\bf f}^0 \|^2= \biggl [ { \Lambda_{\rm sparse}^2 
 \over \phi^4 ( 2, S_0 , 2 s_0) } \biggr ] O ( \lambda_{\rm init}^2 s_0) , $$
 $$ \| b^{S_{\rm init}^{\delta}} - b^0 \|_2=\biggl [  { \Lambda_{\rm sparse} (s_0) \over \phi_{\rm sparse} ( S_0 , 2 s_0) \phi^2 ( 2 , S_0 , 2 s_0) } \biggr ] O ( \lambda_{\rm init} \sqrt {s_0} ) , $$
 and
 $$| S_{\rm init}^{\delta} \backslash S_0 | = O (s_0) . $$
\end{theorem}

\subsection{The adaptive Lasso}\label{adap.section}
Observe that the adaptive Lasso is somewhat more reluctant than thresholding and
refitting: the latter ruthlessly disregards all coefficients with
$| \beta_{j, {\rm init}} | \le \delta $ (i.e., these coefficients get penalty
$\infty$), 
and puts zero penalty on coefficients with
$| \beta_{j, {\rm init}} | > \delta $. The adaptive Lasso gives the coefficients
with $| \beta_{j, {\rm init}} | \le \delta $ a penalty of at least
$\lambda_{\rm init} (\lambda_{\rm adap} /\delta)$ and those
with $| \beta_{j, {\rm init}} | >  \delta $ a penalty of at most
$\lambda_{\rm init} (\lambda_{\rm adap} /\delta)$. 
(Looking ahead, we will actually need to choose $\lambda_{\rm adap} \ge \delta$ in the noisy
case, see Theorem \ref{noisyadap.theorem}.)

Recall
$$\delta_{\rm adap} := \| f_{\rm adap} - {\bf f}^0 \| . $$

The noiseless versions of Conditions B and BB are: 

{\bf Condition b} {\it We have 
$$ \lambda_{\rm init}\biggl [  { \phi_{\rm min} ( 2 , S_0 , 2s_0) \Lambda_{\rm sparse}
 (s_0) \over  \phi^4 ( 2 , S_0 , 2s_0) } \biggr ] = O_{\rm suff} (\lambda_{\rm adap} ) . $$} 
 
  {\bf Condition bb} {\it We have
$$ \lambda_{\rm init}\biggl [ { \phi_{\rm min} ( 2 , S_0 , 2s_0) \Lambda_{\rm sparse}
 (s_0) \over  \phi^4 ( 2 , S_0 , 2s_0) }\biggr ]   \asymp_{\rm suff}  \lambda_{\rm adap}   .$$ }

Note the slight discrepancy with the noisy versions:
the noiseless versions are somewhat better. This is due to
the fact that we also will need to choose $\lambda_{\rm adap}$ large enough
to handle the noise.
 
\begin{theorem}\label{adap.theorem}
Assume Condition b. Then
 $$\delta_{\rm adap}^2 = \biggl [ { \Lambda_{\rm sparse} (s_0) \over
 \phi_{\rm min} ( 2 , S_0 , 2 s_0) }  \biggr ]  { \lambda_{\rm adap}
 \over \lambda_{\rm init} } O ( \lambda_{\rm init}^2 
 s_0 ) , $$
 and
 $$ \|\beta_{\rm adap} - b^0 \|_1 =
 \biggl [ { \Lambda_{\rm sparse}^{1/2} (s_0) \over
 \phi_{\rm min}^{3/2} ( 2 , S_0 , 2 s_0) }  \biggr ] 
 \sqrt { \lambda_{\rm adap}\over \lambda_{\rm init} }O ( \lambda_{\rm init} 
 s_0 ) , $$
 and
  $$ \|\beta_{\rm adap} - b^0 \|_2 =
 \biggl [ { \Lambda_{\rm sparse}^{1/2}  (s_0) 
 \phi_{\rm min}^{1/2} (2 , S_0 ,2 s_0) \over
 \phi_{\rm min}^2 ( 2 , S_0 , 3 s_0) }  \biggr ] \sqrt { \lambda_{\rm adap}
 \over \lambda_{\rm init} } 
 O ( \lambda_{\rm init} 
\sqrt {s_0 } ) , $$
 and
 $$| S_{\rm adap} \backslash S_0|  = 
{  \Lambda_{\rm sparse}^2 (s_0) \over
\phi^4 ( 2 , S_0 , 2 s_0) }  \biggl [  {\Lambda_{\rm sparse} (s_0) \over
\phi_{\rm min} ( 2 , S_0 , 2 s_0) } \biggr ] {\lambda_{\rm init} 
\over \lambda_{\rm adap} } O ( s_0) .  $$
Considering the bounds for the
prediction and estimation error  leads to favoring the choice of Condition bb, giving
  $$\delta_{\rm adap}^2 = \biggl [ { \Lambda_{\rm sparse}^2 (s_0) \over
 \phi^4 ( 2 , S_0 , 2 s_0) }  \biggr ] O ( \lambda_{\rm init}^2 
s_0 ) , $$
 $$ \|\beta_{\rm adap} - b^0 \|_1 =
 \biggl [ { \Lambda_{\rm sparse} (s_0) \over
 \phi_{\rm min} ( 2 , S_0 , 2 s_0) \phi^2(2, S_0 , 2 s_0)}  \biggr ] O (\lambda_{\rm init} 
 s_0 ) , $$
 $$ \|\beta_{\rm adap} -b^0 \|_2 =
 \biggl [ { \Lambda_{\rm sparse} (s_0) \phi_{\rm min}( 2 , S_0 , 2 s_0) 
 \over
 \phi_{\rm min}^2 ( 2 , S_0 , 3 s_0) \phi^2 (2 , S_0 , 2 s_0)  }  \biggr ] O ( \lambda_{\rm init} 
\sqrt {s_0 } ) , $$
 and 
  $$| S_{\rm adap} \backslash S_0|  = 
{  \Lambda_{\rm sparse}^2 (s_0) \over
\phi_{\min}^2 ( 2 , S_0 , 2 s_0) }  O ( s_0) .  $$
\end{theorem} 
 
\subsection{The weighted irrepresentable condition}\label{irrepresentable.section}

This subsection will show that, even in the noiseless case, exact variable selection
 needs rather strong conditions.
It serves as a motivation for the perhaps more moderate aim of having
$O(s_0) $ ($\le O(s_{\rm true})$) false positives and detecting only the larger
coefficients. Moreover, we illustrate in Example \ref{separate.example}
of this subsection
that the lower bound on the non-zero coefficients
as given in Corollary \ref{min.corollary} is tight.

It is known that the initial Lasso essentially needs the irrepresentable
condition in order to have no false positives  ({\cite{ZY07}).
Similar statements
can be made for the weighted Lasso.

 For a $(p \times p)$ -matrix $\Sigma= (\sigma_{j,k})$. we define
 $$\Sigma_{1,1} (S) := (\sigma_{j,k})_{j,k \in S}, $$
 $$\Sigma_{2.1} (S):= (\sigma_{j,k})_{j \notin S ,\ k \in S} . $$
 We let $W_S := {\rm diag } ( \{ w_j \}_{j \in S} ) $.

 {\bf Definition} {\it We say that the {\rm weighted irrepresentable
 condition} holds for $S$ if for all vectors $\tau_{S} \in \R^{|S|}$ with
 $\| \tau_{S} \|_{\infty} \le 1$, one has
 $$\| W_{S^c}^{-1} \Sigma_{2,1} (S) \Sigma_{1,1}^{-1} (S) W_{S} \tau_{S} \|_{\infty}< 1 . $$
 }

The reparametrization $\beta \mapsto
\gamma := W^{-1} \beta$ leads to the following lemma,
which is the weighted variant of the first part of Lemma 6.2 in
\cite{vdG:2009}. Here, we actually take ${\bf f}_0$ as target,
instead of its $\ell_0$-sparse approximation ${\rm f}_{S_0}$.
Recall
$$S_{\rm true} := \{ j : \ \beta_{j, {\rm true} } \not= 0 \} . $$

\begin{lemma}\label{irrepresentable}$ $\\
Suppose the weighted irrepresentable condition is met for $S_{\rm true}$. Then $S_{\rm weight} \subset S_{\rm true}$.
\end{lemma}

We now consider conditions for the weighted irrepresentable
condition to hold.

\begin{lemma}\label{boundirr.lemma}
Suppose that
\begin{equation} \label{boundirr}
 \| w_S \|_2 < \Lambda_{\rm min}(S) w_{S^c}^{\rm min} . 
 \end{equation}
Then the weighted irrepresentable condition holds for $S$.
\end{lemma}

The next example shows that the result of Lemma \ref{boundirr.lemma}
cannot be improved (essentially, up to the strict inequality)
without assuming further conditions.

\begin{example} \label{separate.example} {
Let $S_{\rm true}= \{ 1 , \ldots , s \}$, with cardinality $s := |S_{\rm true}|$, be the active set, and write
 $$ \Sigma := \pmatrix { \Sigma_{1,1} &  \Sigma_{1,2} \cr  \Sigma_{2,1} & \Sigma_{2,2} \cr}.
$$
We now will take a special choice for $\Sigma$, which is
perhaps not very representative when $\Sigma$ is an empirical
Gram matrix $\hat \Sigma$, but it is legitimate for a worst case analysis
(as we study here). 
We suppose that  $ \Sigma_{1,1}:=I$ is the $(s \times s)$-identity matrix,
and
$$ \Sigma_{2,1}:= \rho( c_2 c_1^T)  , $$
with $0\le  \rho < 1$, and with $c_1$ an $s$-vector and $c_2$ a $(p-s)$-vector,
satisfying $\| c_1 \|_2=\| c_2 \|_2=1 $.
Moreover, we suppose
$ \Sigma_{2,2} $ is the $((p-s)\times (p-s) )$-identity matrix.
Then $\Lambda_{\rm min} (S_{\rm true} )=1$ and the smallest eigenvalue of $\Sigma$ is
$1- \rho$. Its largest eigenvalue is $1+ \rho$.
Take $c_1= w_{S_{\rm true}}/ \| w_{S_{\rm true}} \|_2$,
 and $c_2 = (0 , \ldots , 1 , 0 , \ldots )^T$,
 where the $1$ is placed at $\arg \min_{j \in S_{\rm true}^c} w_j $.
 Then
 $$
 \sup_{\| \tau_{S_{\rm true}} \|_{\infty} \le 1 }
 \| W_{S_{\rm true}^c}^{-1}\Sigma_{2,1} \Sigma_{1,1}^{-1}W_{S_{\rm true}}  \tau_{S_{\rm true}} \|_{\infty} =
 \rho \| w_{S_{\rm true}} \|_2 / w_{S_{\rm true}^c}^{\rm min}  .$$
 As a special case, suppose $c_1 = (1 , 1 , \ldots , 1 )^T / \sqrt {s}$, and 
 $\rho=1/2$.
The adaptive Lasso generally has $O(1/ w_{S_{\rm true}^c}^{\min})
 = \lambda_{\rm init}$. The irrepresentable
 condition then needs
 $$\sum_{j \in S_{\rm true} } 1/ |\beta_{j, {\rm init}} |^2 =
 \| w_{S_{\rm true}} \|_2^2 = O (\lambda_{\rm init}^{-2} )$$ 
 (which holds for example
 when $ \lambda_{\rm init} \sqrt {s} = O(\min_{j \in S_{\rm true}} | \beta_{j, {\rm true}} |) $.)
 This condition also shows up in Corollary \ref{min.corollary},
 i.e., the result there is tight.}
\end{example}
 
 \section{Adding noise} \label{noise.section}
 
 After introducing the notation for the noisy case (Subsection
 \ref{noisynotation.section}),
 we will give the extension of the results for the weighted Lasso
 to the noisy case\footnote{Of separate interest is a direct comparison of the noisy initial Lasso with
the noisy $\ell_0$-penalized estimator. Replacing ${\bf f}^0$ by ${\bf Y}$ in Corollary \ref{init.corollary} (and dropping the requirement $S \subset S_{\rm true}$)
gives
$$\| {\bf Y} - \hat f_{\rm init} \|_n^2 \le 2 \min_{S} \left \{ 
\| {\bf Y} - \hat {\rm f}_S \|_n^2 + {3 \lambda_{\rm init}^2 | S | \over
\phi^2 ( 2, S) } \right \} . $$}
  (see Theorem \ref{weightnoise}).
Once this is done, results for the initial Lasso, its
thresholded version, and for the adaptive Lasso, follow in the
same way as in Subsections \ref{init.section}, \ref{thres.section} and  \ref{adap.section}.
The new point is to take care that the tuning
parameters are chosen in such a way that the
noisy part due to variables in $S_0^c$ are overruled
by the penalty term. In our situation, this can be done
by taking $\lambda_{\rm init}$, as well as
$\lambda_{\rm adap}\ge \lambda_{\rm init}$ sufficiently large.

We provide the result for the noisy weighted Lasso in Subsection
\ref{noisyweight.section}.  Theorems \ref{noisyinit.theorem},
\ref{noisythres.theorem} and \ref{noisyadap.theorem} follow from this
and from some further results for the noisy case (their proofs are in Subsection
\ref{preview.proofs}).
In Section \ref{refine.section}, we look at more restrictive sparse
eigenvalue conditions in the spirit of
\cite{zhang2008sparsity}.

\subsection{Notation for the noisy case}\label{noisynotation.section}
Consider an $n$-dimensional vector of observations 
$$
{\bf Y} = {\bf f}^0  + \epsilon.
$$
where ${\bf f}^0 := ({\bf f}^0 (X_1 ) , \ldots , {\bf f}^0 (X_n))^T $,
with $X_1 , \ldots , X_n$ co-variables in some space ${\cal X}$.
Let $\{ \psi_j \}_{j=1}^p $ be a given dictionary.

The regression ${\bf f}^0$, the dictionary $\{ \psi_j \}$, and $f_{\beta} := \sum \psi_j \beta_j $ are now
considered as vectors in $\R^n$.
The norm we use is
the normalized Euclidean norm
$$\| f \| := \| f \|_n:= \| f \|_2/ \sqrt {n } : \ f \in {\bf R}^n , $$
induced by the inner product
$$( f , \tilde f)_n :=  {1 \over n} \sum_{i=1}^n f_i \tilde f_i  , \ 
f , \tilde f \in \R^n  . $$
In other words, the probability measure $Q$ is now
$Q:= Q_n = \sum_{i=1}^n \delta_{X_i}/n$, the empirical
measure of the co-variables $X_1 , \ldots , X_n$.
With some abuse of notation, we also write
$$\| {\bf Y} - f \|_n^2 := \| {\bf Y} - f \|_2^2 / n , $$
and
$$(\epsilon , f)_n := {1 \over n} \sum_{i=1}^n  \epsilon_i f (X_i) . $$

The design matrix ${\bf X} $ is 
$${\bf X} = ( \psi_1 , \ldots , \psi_p ) . $$
We write the eigenvalues
involved as before, e.g.,
$\Lambda_{\rm max} $ is the largest eigenvalue of the 
empirical Gram matrix $\hat \Sigma := {\bf X}^T {\bf X} / n $, and
$ \phi^2 (L,S,N)$ is the $(L,S,N)$-restricted eigenvalue of $\hat \Sigma$.
The projections in $L_2 (Q_n)$ are also written as before, i.e.
$${\rm f}_S: = {\bf X} b^S := \arg \min_{f = {\bf X} \beta_S } \| f - {\bf f}^0 \|_n . $$

 The $\ell_0$-sparse projection ${\rm f}_{S_0}= \sum_{j \in S_0}
b_j^0 $ is now defined with a larger constant  (7 instead of 3) in front of the
penalty term, and a larger constant ($ L=6 $ instead of $L=2$) in the restrictions
 of the restricted eigenvalue condition:
$$S_0 := \arg \min_{S \subset S_{\rm true}} \left \{
\| {\rm f}_S - {\bf f}^0 \|_n^2 + { 7 \lambda_{\rm init}^2 |S| \over
 \phi^2 ( 6, S) } \right \}  $$
(compare with formula (\ref{defineS_0})).

The weighted Lasso is
\begin{equation}
\label{eq.weight}
\hat \beta_{\rm weight} = \arg\min_{\beta} \biggl \{  \|{\bf Y} -f_{\beta} \|_n^2 +  
\lambda_{\rm init} \lambda_{\rm weight}  \sum_{j=1}^p w_j |\beta_j |
\biggr \} .
\end{equation}
Let
$$\hat f_{\rm weight} := f_{\hat \beta_{\rm weight}}, \ 
\hat S_{\rm weight} := \{ j: \ \hat \beta_{j, {\rm weight}} \not= 0 \} . $$

The initial and adaptive Lasso are defined as in Section \ref{introduction.section}.
We write $\hat f_{\rm init} := f_{\hat \beta_{\rm init}}$ and $\hat f_{\rm adap}
:= f_{\hat \beta_{\rm adap}}$, with active sets $\hat S_{\rm init}:=
\{j:\  \hat \beta_{j,\rm init} \not= 0 \}$ and
$\hat S_{\rm adap}:= \{ j:\ \hat \beta_{j,\rm adap} \not= 0 \} $, respectively. 
Let
$$ \hat \delta_{\rm init}^2 := \| f_{\beta_{\rm init}} - {\bf f}^0  \|_n^2 , $$
be the prediction error of the initial Lasso, and 
and, for $q \ge 1 $,
$$\hat \delta_q := \| \hat \beta_{\rm init} - b^0 \|_q  $$
be its $\ell_q$-error.
Denote the prediction error of the adaptive Lasso by
$$ \hat \delta_{\rm adap}^2 :=  \| f_{ \hat \beta_{\rm adap}} - {\bf f}^0  \|_n^2 .$$

The least squares estimator using only variables in $S$ is also written with a ``hat":
$$\hat {\rm f}_S  =
f_{\hat b^S } := \arg \min_{f = f_{\beta_S}}  \| {\bf Y} - f_{\beta_S } \|_n . $$
A threshold level will be denoted by $\delta$, instead of $\lambda_{\rm thres}$
as we do in Section \ref{introduction.section}.
The reason is again that we need to explicitly express dependence
on the threshold level. With $\lambda_{\rm thres}$ the notation will
be too complicated.
We define, for any threshold $\delta > 0$, 
$$\hat S_{\rm init}^{\delta} := \{ j : \ | \hat \beta_{j,{\rm init}} | > \delta \} .
  $$
The refitted version after thresholding, 
based on the data ${\bf Y}$, is
$ \hat {\rm f}_{\hat S_{\rm init}^{\delta}} $.

To handle the (random) noise, we define the set
$${\cal T} := \biggl  \{ \max_{1 \le j \le p} 4 | ( \epsilon , \psi_j )_n | \le 
\lambda_{\rm init} \biggr \}  .$$
This is the set where the (empirical) correlations between noise and
design is ``small".

Here $\lambda_{\rm init}$ is chosen in such a way that
$$\PP ( {\cal T} ) \ge 1- \alpha $$
where $(1-\alpha)$ is the confidence we want to achieve.

\subsection{The noisy weighted Lasso} \label{noisyweight.section}

\begin{theorem} \label{weightnoise}
Suppose we are on ${\cal T}$. 
Let $S$ be a set with cardinality $s= |S|$, which satisfies for some 
positive $L$ and $M$
$$
\lambda_{\rm weight} ( w_{S^c}^{\rm min}   \wedge M)  \ge 1,$$
and
$$ w_{S^c}^{\rm min}  \ge M/L, \ \| w_S \|_2/\sqrt {s} \le M  . $$
Then
for all $\beta$,
$$\| \hat  f_{\rm weight} - {\bf f}^0 \|_n^2    \le
2 \| f_{\beta_S} - {\bf f}^0 \|_n^2 + { 14 \lambda_{\rm init}^2
\lambda_{\rm weight}^2  M^2  s \over \phi^2 (6L ,S) } 
 ,$$
and
$$  \sqrt {s} 
\| (\hat \beta_{\rm weight})_S - \beta_S \|_2 +  
 \| (\hat  \beta_{ {\rm weight}})_{S^c} \|_1/L $$
$$ \le 
{ 5 \| f_{\beta_S} - {\bf f}^0 \|_n^2
\over \lambda_{\rm init} \lambda_{\rm weight} M}   +{ 7 \lambda_{\rm init}
\lambda_{\rm weight}  M s \over \phi^2 (6L,S) }  ,$$
and
$$\| \hat \beta_{\rm weight} - \beta_S \|_2 $$ $$ \le
 { 10L \| f_{\beta_S}  - {\bf f}^0 \|_n^2  \over M  \lambda_{\rm init}
\lambda_{\rm weight} \sqrt {s_0} } 
 + {14 L\lambda_{\rm init}^2 \lambda_{\rm weight}^2M  (s  +s_0) \over   \phi^2 (6L, S, s+s_0) \lambda_{\rm init}
\lambda_{\rm weight}   \sqrt {s_0}}  . 
 $$
 Moreover,
 under the condition $\lambda_{\rm weight} w_{S^c}^{\rm min} \ge 1$, 
$$|( \hat S_{\rm weight} \cap S^c)  \backslash S_0 |^2 $$ $$\le 16
 \Lambda_{\rm max}^2 (( \hat S_{\rm weight} \cap S^c)  \backslash S_0)
 {\| \hat f_{\rm weight} - {\bf f}^0 \|_n^2   \over \lambda_{\rm weight}^2  }
 {  \|( 1/w)_{\hat S_{\rm weight} \backslash S_0}\|_2^2 \over \lambda_{\rm init}^2}. $$
When $|( \hat S_{\rm weight} \cap S^c)  \backslash S_0 |> s_0$,
this implies
$$| (\hat S_{\rm weight} \cap S^c ) \backslash S_0 | \le 
32 \Lambda_{\rm sparse}^2 (s_0)  {\| \hat f_{\rm weight} - {\bf f}^0 \|_n^2   \over \lambda_{\rm weight}^2 s_0 }
{  \|( 1/w)_{\hat S_{\rm weight} \backslash S_0}\|_2^2 \over \lambda_{\rm init}^2}
.$$

\end{theorem}

\subsection{Another look at the number of false positives} \label{refine.section}

Here, we discuss a refinement, assuming a condition
corresponding to the one used in \cite{zhang2008sparsity}.

{\bf Condition D} {\it It holds for some $s_* \ge s_0$, that
$$D (s_*, s_0) := \biggl \{ { \Lambda_{\rm sparse}^2 (s_* ) s_0 \over 
\phi^2 ( 2, S_0) s_*} \biggr \} = O_{\rm suff} (1) . $$} 

\begin{lemma} \label{refineselect.lemma} Suppose we are on ${\cal T}$.
Then under Condition D,
$$ |\hat S_{\rm init}  \backslash S_0 | = \biggl [ 
{ \Lambda_{\rm sparse}^2 (s_*) \over
\phi^2 ( 6 , S_0) } \biggr ] \biggl ( 1- {D(s_* , s_0)
\over O_{\rm suff} (1)}  \biggr )^{-1}  O(s_0) . $$
Moreover, under Condition B, 
$$|\hat S_{\rm adap} \backslash S_0 | =
  \Lambda_{\rm sparse} ( s_*) \biggl [ { \Lambda_{\rm sparse} (s_0) \over \phi_{\rm min} ( 6 , S_0 , 2 s_0)  
\phi^4 ( 6 , S_0 , 2 s_0) }  \biggr ]^{1/2} 
\sqrt { \lambda_{\rm init} \over \lambda_{\rm adap} }
O(s_0)   $$
$$ +
\biggl [ { \Lambda_{\rm sparse} (s_0)  \phi^2 ( 6, S_0 ) \over \phi_{\rm min} ( 6 , S_0 , 2 s_0)  
\phi^4 ( 6 , S_0 , 2 s_0)}  \biggr ] 
D(s, s_*)  {\lambda_{\rm init} 
\over \lambda_{\rm adap} } O(s_0) .$$ 
Under Condition BB, this becomes
\begin{equation} \label{adap**}
|\hat S_{\rm adap} \backslash S_0 | =
 \biggl [  { \Lambda_{\rm sparse} ( s_*) \over 
\phi(6, S_0) }\biggr ]  \biggl [ { \phi_{\rm min}^2 ( 6 , S_0 , 2 s_0) 
\phi^2 (6, S_0) \over \phi^2 ( 6 , S_0 , 2 s_0)  
 }  \biggr ]^{1/2} 
O(s_0 ) 
\end{equation}
$$ + 
\biggl [ { \phi_{\rm min}^2 ( 6 , S_0 , 2 s_0)  \phi^2 ( 2, S_0 ) \over
\phi^4 ( 6 , S_0 , 2 s_0)}  \biggr ]
D(s_*, s_0)   O(s_0). $$
\end{lemma}

Under Condition D, the first term in the right hand side of 
(\ref{adap**}) is generally the leading term. We thus see the adaptive Lasso
replaces the potentially very large constant 
$$ \biggl ( 1- {D(s_* , s_0)
\over O_{\rm suff} (1)}  \biggr )^{-1} $$in the
bound for the number of false positives of the initial Lasso by 
$$
\biggl [ { \phi_{\rm min}^2 ( 6 , S_0 , 2 s_0) 
\phi^2 (6, S_0) \over \phi^4 ( 6 , S_0 , 2 s_0)  
 }  \biggr ]^{1/2} , $$
 a constant which is close to 1 if the $\phi$'s do not differ too much.
 
 Admittedly, Condition D is difficult to interpret. On the one hand, it 
 wants $s_*$ to be large, but on the other hand, a large $s_*$ also
 can render $\Lambda_{\rm sparse} (s_*)$ large. 
We refer to \cite{zhang2008sparsity} for examples where
Condition D is met.

\section{Proofs}\label{proofs.section}

We present three subsections, containing respectively
the proofs for Section \ref{noiseless.section}, Section \ref{noise.section},
and finally Section \ref{preview.section}. 

\subsection{Proofs for Section \ref{noiseless.section}: the noiseless case}

\subsubsection{Proofs for Subsection \ref{weight.section}: the noiseless weighted Lasso}
{\bf Proof of Theorem \ref{weight}.} 
Take
$$   w_{S^c}^{\rm min}  \ge M/L,  \  \| w_S \|_2 / \sqrt {s} \le M . $$
We have
$$\|  f_{\rm weight} - {\bf f}^0 \|^2 +\lambda_{\rm init}
\lambda_{\rm weight}  \sum_{j=1}^p w_j | \beta_{j, {\rm weight}}  | \le \| f_{\beta_S} - {\bf f}^0\|^2 +
\lambda_{\rm init}
\lambda_{\rm weight} \sum_{j \in S } w_j |\beta_j | , $$
and hence
$$\|  f_{\rm weight} - {\bf f}^0 \|^2 + \lambda_{\rm init}
\lambda_{\rm weight} w_{S^c}^{\rm min}  \| ( \beta_{ {\rm weight}})_{S^c} \|_1 
$$ $$\le
 \| f_{\beta_S}- {\bf f}^0 \|^2 +
\lambda_{\rm init}
\lambda_{\rm weight}  \sum_{j \in S} w_j |  \beta_{j, {\rm weight}}- \beta_j |  $$
$$\le
 \| f_{\beta_S}- {\bf f}^0 \|^2 +
\lambda_{\rm init}
\lambda_{\rm weight}  M \sqrt {s}  \| ( \beta_{\rm weight})_S  - \beta_S \|_2 . $$
Let ${\cal N} \supset S$, $|{\cal N} | = N$. Then
$$\|( \beta_{\rm weight})_{{\cal N}^c} \|_1 \le \| (\beta_{\rm weight})_{S^c} \|_1 , $$
and $$
\| ( \beta_{\rm weight})_S - \beta_S \|_2 \le 
\| ( \beta_{\rm weight})_{\cal N}  - \beta_S \|_2 , \ \sqrt {s} \le \sqrt N .$$
Therefore, 
$$\|  f_{\rm weight} - {\bf f}^0 \|^2 + \lambda_{\rm init}
\lambda_{\rm weight} w_{S^c}^{\rm min}  \| ( \beta_{ {\rm weight}})_{{\cal N}^c} \|_1 $$
$$\le
 \| f_{\beta_S}- {\bf f}^0 \|^2 +
\lambda_{\rm init}
\lambda_{\rm weight}  M \sqrt {N}  \| ( \beta_{\rm weight})_{\cal N}  - \beta_S \|_2 . $$

{\bf Case i).} If
$$  \| f_{\beta_S}- {\bf f}^0 \|^2 \le 
\lambda_{\rm init}
\lambda_{\rm weight} M \sqrt {N} \| ( \beta_{\rm weight})_{\cal N}  - \beta_S \|_2, $$
we get
\begin{equation}\label{important}
\|  f_{\rm weight} - {\bf f} ^0 \|^2 + \lambda_{\rm init}
\lambda_{\rm weight} w_{S^c}^{\rm min}  \| ( \beta_{ {\rm weight}})_{{\cal N}^c} \|_1 
\end{equation}
$$
\le 2
\lambda_{\rm init}
\lambda_{\rm weight}  M \sqrt {N}  \| ( \beta_{\rm weight})_{\cal N}  - \beta_S \|_2. $$

It follows that
$$\|( \beta_{\rm weight})_{{\cal N}^c} \|_1 \le  2L\sqrt {N} 
\|(\beta_{\rm weight})_{\cal N} - (\beta)_{S} \|_2. $$
But then, by the definition of restricted eigenvalue, and invoking the
triangle inequality,
$$  \|  (\beta_{\rm weight})_{\cal N}- \beta_S \|_2 
 \le  \| f_{\rm weight} - f_{\beta_S}  \| / \phi (2L,{\cal N})  $$
$$ \le  \| f_{\rm weight} - {\bf f}^0 \| / \phi (2L,{\cal N}) +
 \| f_{\beta_S} - {\bf f}^0 \| / \phi (2L,{\cal N}) .$$
This gives
$$\|  f_{\rm weight} - {\bf f} ^0 \|^2 + \lambda_{\rm init}
\lambda_{\rm weight} w_{S^c}^{\rm min}  \| ( \beta_{ {\rm weight}})_{{\cal N}^c} \|_1 
 $$ $$\le 2
\lambda_{\rm init}
\lambda_{\rm weight} 
 M \sqrt {N}\|  f_{\rm weight} - {\bf f}^0 \| / \phi (2L,{\cal N})  $$ $$+
 2\lambda_{\rm init}
\lambda_{\rm weight} 
M\sqrt N \| f_{\beta_S} - {\bf f}^0 \| / \phi (2L,{\cal N}) $$
 $$\le {1 \over 2}  \| f_{\rm weight} - {\bf f}^0 \|^2 + 
 \| f_{\beta_S} - {\bf f}^0 \|^2 + { 3 \lambda_{\rm init}^2
\lambda_{\rm weight}^2  N M^2\over
 \phi^2 (2L,{\cal N} ) } . $$
 Hence,
 $$
 \| f_{\rm weight}  - {\bf f} ^0 \|^2 + 2 \lambda_{\rm init}
\lambda_{\rm weight}  w_{S^c}^{\rm min} 
\| (\beta_{\rm weight})_{{\cal N}^c} \|_1\le 2 \| f_{\beta_S} - {\bf f}^0 \|^2 + { 6 \lambda_{\rm init}^2
\lambda_{\rm weight}^2     N M^2\over
 \phi^2 (2L,{\cal N} ) } . $$
 
 {\bf Case ii)} 
 If
$$  \| f_{\beta_S}- {\bf f}^0 \|^2 >
\lambda_{\rm init}
\lambda_{\rm weight}  M \sqrt N   \| ( \beta_{\rm weight})_{\cal N}  - \beta_S \|_2
 , $$
we get
$$\|  f_{\rm weight}- {\bf f} ^0 \|^2 + \lambda_{\rm init}
\lambda_{\rm weight} w_{S^c}^{\rm min}  \| ( \beta_{ {\rm weight}})_{{\cal N}^c} \|_1 
\le 
2 \| f_{\beta_S} - {\bf f}^0 \|^2 . $$

The first result of the Lemma now follows from taking ${\cal N} = S$.

For the second result, we add in Case i), 
$\lambda_{\rm init}
\lambda_{\rm weight}  M  \sqrt N \| ( \beta_{\rm weight})_{\cal N}  - \beta_S \|_2$
to the left and right hand side of (\ref{important}):
$$\|  f_{\rm weight} - {\bf f} ^0 \|^2 + \lambda_{\rm init}
\lambda_{\rm weight} M \sqrt N \| ( \beta_{\rm weight})_{\cal N}  - \beta_S \|_2 $$ $$+
 \lambda_{\rm init}
\lambda_{\rm weight} w_{S^c}^{\rm min}  \| ( \beta_{ {\rm weight}})_{{\cal N}^c} \|_1  
 $$ $$
\le 3
\lambda_{\rm init}
\lambda_{\rm weight}  M \sqrt N \| ( \beta_{\rm weight})_{\cal N}  - \beta_S \|_2. $$
The same arguments now give
$$ 3
\lambda_{\rm init}\sqrt N\| ( \beta_{\rm weight})_{\cal N}  - \beta_S \|_2
+
 \lambda_{\rm init}
\lambda_{\rm weight} w_{S^c}^{\rm min}  \| ( \beta_{ {\rm weight}})_{{\cal N}^c} \|_1  \le  $$
$$ \| f_{\rm weight} - {\bf f}^0 \|^2 +
3 \| f_{\beta_S} - {\bf f}^0 \|^2 + { 3 \lambda_{\rm init}^2 \lambda_{\rm weight}^2 
N M^2   \over
\phi^2 ( 2 L, {\cal N}) } . $$

In Case ii), we have
$$ \lambda_{\rm init}
\lambda_{\rm weight} w_{S^c}^{\rm min}  \| ( \beta_{ {\rm weight}})_{{\cal N}^c} \|_1  
\le 
2 \| f_{\beta_S} - {\bf f}^0 \|^2  , $$
and also
$$\lambda_{\rm init}
\lambda_{\rm weight}  M \sqrt N  \| ( \beta_{\rm weight})_{\cal N}  - \beta_S \|_2 <
 \| f_{\beta_S}- {\bf f}^0 \|^2 . $$
 So then
 $$
 \lambda_{\rm init}
\lambda_{\rm weight}  M \sqrt N  \| ( \beta_{\rm weight})_{\cal N}  - \beta_S \|_2 
 \lambda_{\rm init}
\lambda_{\rm weight} w_{S^c}^{\rm min}  \| ( \beta_{ {\rm weight}})_{{\cal N}^c} \|_1  +
$$ $$<
3 \| f_{\beta_S} - {\bf f}^0 \|^2 . $$
Taking ${\cal N} = S$ gives the second result.

For the third result, 
we let 
${\cal N}$ be the set $S$, complemented with the $s_0$ largest - in absolute value - coefficients
of $(\beta_{\rm weight})_{S^c}$. Then
$\phi(2L, {\cal N}) \le \phi( 2, S , s+ s_0). $
Moreover, $N \ge  s_0$.
Thus, from the second result, we get
$$\lambda_{\rm init} \lambda_{\rm weight} M \sqrt { s_0} \| (\beta_{\rm weight})_{\cal N} -\beta_S \|_2 
+ \lambda_{\rm init} \lambda_{\rm weight} \| ( \beta_{\rm weight} )_{{\cal N}^c} \|_1 $$ $$ \le
3 \| f_{\beta_S} - {\bf f}^0 \|^2 + {3 \lambda_{\rm init}^2 \lambda_{\rm weight}^2 (s_0+s) M^2 \over
\phi^2 (2L, S, s+s_0 ) } . $$
Moreover, as is shown in Lemma 2.2 in  \cite{vdG:2009} (with original reference
\cite {candes2005decoding}, and \cite{candes2007dss}), 
$$\| (\beta_{\rm weight})_{{\cal N}^c} \|_2 \le
\| (\beta_{\rm weight})_{S^c}  \|_1 / \sqrt {s_0}  $$ $$
\le { 3L \| f_{\beta_S} - {\bf f}^0 \|^2 + 3L \lambda_{\rm init}^2 (s+s_0) M^2 / \phi^2 (L, S , s+s_0)  \over 
 \lambda_{\rm init} \lambda_{\rm weight}  M \sqrt {s_0}} .$$
So then
$$\| \beta_{\rm weight} - \beta_S \|_2 \le
\| (\beta_{\rm weight})_{\cal N} - \beta_S \|_2 +
\| ( \beta_{\rm weight})_{{\cal N}^c }  \|_2 $$
$$ \le { 6L \| f_{\beta_S}  - {\bf f}^0 \|^2 + 6L \lambda_{\rm init}^2 \lambda_{\rm weight}^2 (s+ s_0)M^2  / \phi^2 (2L, S, s+s_0)  \over M \sqrt {s_0} \lambda_{\rm init}
\lambda_{\rm weight} }  . 
 $$
\hfill $\sqcup \mkern -12mu \sqcap$

We now turn to the proof of Lemma \ref{weightselect}. 
An important characterization of the solution $\beta_{\rm weight}$ can be derived from
 the {\it Karush-Kuhn-Tucker} ({\it KKT}) conditions (see \cite{bertsimas1997introduction}).
 
{\bf Weighted KKT-conditions}
{\it We have 
$$2 \Sigma ( \beta_{\rm weight} - \beta_{\rm true} ) =-\lambda_{\rm weight} \lambda_{\rm init} 
W \tau_{\rm weight} .$$
Here, $\| \tau_{\rm weight} \|_{\infty} \le 1$, and moreover
$$\tau_{j, {\rm weight}} {\rm l} \{ \beta_{j, {\rm weight}} \not= 0 \}= {\rm sign} (\beta_{j, {\rm weight}} ) , \ j = 1 , \ldots , p. $$}

{\bf Proof of Lemma \ref{weightselect}.}
By the weighted KKT conditions, for all $j$
$$2 ( \psi_j , f_{\rm weight} - {\bf f}^0 )  =- \lambda_{\rm init} \lambda_{\rm weight} w_j \tau_{j, {\rm weight}} .$$
Hence, 
$$\sum_{j \in S_{\rm weight} \backslash S_0 }
2 | ( \psi_j , f_{\rm weight}  - {\bf f}^0 ) |^2 \ge \lambda_{\rm init}^2 
\lambda_{\rm weight}^2\| w_{S_{\rm weight} \backslash S_0 } \|_2^2 $$ $$
\ge\lambda_{\rm init}^2 
\lambda_{\rm weight}^2 | S_{\rm weight} \backslash S_0 |^2/  \|( 1/w)_{S_{\rm weight} \backslash S_0}\|_2^2
. $$
On the other hand
$$\sum_{j \in S_{\rm weight} \backslash S_0 }
 | ( \psi_j , f_{\rm weight} - {\bf f}^0 ) |^2 \le \Lambda_{\rm max }^2 (S_{\rm weight} \backslash
S_0) \| f_{\rm weight} - {\bf f}^0 \| ^2 . $$
Thus, we arrive at inequality (\ref{stillin}):
$$
| S_{\rm weight} \backslash S_0 |^2 \le 4
\Lambda_{\rm max}^2 ( S_{\rm weight} \backslash S_0 ) 
{ \| f_{\rm weight} - {\bf f}^0 \|^2 \over \lambda_{\rm weight}^2  }
{ \| 1 / w_{S_{\rm weight}\backslash S_0 } \|^2 \over \lambda_{\rm init}^2}  .
$$
Clearly,
$$\Lambda_{\rm max}^2  (S_{\rm weight} \backslash S^0) \le
\Lambda_{\rm max}^2 \wedge \biggl ( { | S_{\rm weight} \backslash S_0| \over s_0 } +1
\biggr ) 
\Lambda_{\rm sparse}^2 (s_0) . $$

\hfill $\sqcup \mkern -12mu \sqcap$

\subsubsection{Proofs for Subsection \ref{init.section}: the noiseless initial Lasso}

We first present the corollaries of Theorem \ref{weight}
and Lemma \ref{weightselect} when we apply them to the case
where all the weights are equal to one.

\begin{corollary}\label{init.corollary}
For the initial Lasso, $w_j=1$ for all $j$, so we can apply Corollary \ref{ordered.corollary}
with $\delta=1$ and $S_{\rm weight}^{\delta} = S_0$.
Let
$$\delta_{\rm oracle}^2 := 
 \| {\rm f}_{S_0} - {\bf f}^0 \|^2 + { 3 \lambda_{\rm init}^2 |S_0| \over \phi^2(2, S_0)  } . $$
We have
$$ \delta_{\rm init}^2  \le
2 \| {\rm f}_{S_0} - {\bf f}^0 \|^2 + { 6 \lambda_{\rm init}^2 |S_0| \over \phi^2 (2 , S_0) }
=
2\delta_{\rm oracle}^2.
$$
The estimation error can be bounded as follows:
$$ \delta_1 \le 3 \| {\rm f}_{S_0} - {\bf f}^0 \|^2/ \lambda_{\rm init}  + { 3 \lambda_{\rm init} |S_0| \over \phi^2 (2 , S_0) } \le 
3 \delta_{\rm oracle}^2 / \lambda_{\rm init} , $$
and
$$ \delta_2 \le  \biggl [ { \phi^2 (2, S_0) \over \phi^2 (2 , S_0 , 2 s_0) }\biggr ] { 6 \delta_{\rm oracle}^2 \over  \lambda_{\rm init} \sqrt {s_0}}
   . $$
   Moreover, application of Lemma \ref{weightselect} bounds the number of
   false positives:
$$| S_{\rm init} \backslash S_0 | \le 4 \Lambda_{\rm max}^2 (S_{\rm init}
\backslash S_0 ) 
 { \delta_{\rm init}^2 \over \lambda_{\rm init}^2} .$$
\end{corollary}

{\bf Proof of Theorem \ref{init.theorem}.}\\ 
This is now a direct consequence of Corollary \ref{init.corollary}.
\hfill $\sqcup \mkern -12mu \sqcap$

\subsubsection{Proofs for Subsection \ref{thres.section}: the noiseless thresholded Lasso}
We first provide some explicit bounds.

\begin{lemma} \label{betterthres}
We have
$$ \| ( \beta_{\rm init})_{S_{\rm init}^{\delta} } -
b^0 \|_1 \le  2 \delta_1 + \delta s_0 , $$
and
$$  \| ( \beta_{\rm init})_{S_{\rm init}^{\delta} } -
b^0 \|_2 \le  2 \delta_2 +  \delta \sqrt {s_0 },$$
and
$$ 
\| {\rm f}_{S_{\rm init}^{\delta} }  - {\bf f}^0 \| \le
\| f_{(\beta_{\rm init})_{S_{\rm init}^{\delta}}} - {\bf f}^0 \| $$ $$ \le 
 \| {\rm f}_{S_0} - {\bf f}^0 \| + \sqrt {\biggl  \lceil {\delta_2^2 \over \delta^2 s_0}
 + 1 \biggr \rceil }   \Lambda_{\rm sparse}  (s_0) ( 2 \delta_2 + \delta  \sqrt{ s_0}), 
  $$
  and, for $\delta \ge \delta_2 / \sqrt {s}$, 
  $$\| b^{S_{\rm init}^{\delta}} - b^0 \|_2 \le { \| {\rm f}_{S_{\rm init}^{\delta}} - {\bf f}^0 \| 
  \over \phi_{\rm sparse} ( S_0,  2s_0) } . $$
\end{lemma}

{\bf Proof of Lemma \ref{betterthres}.}
To obtain the first result, we use
$$  \| ( \beta_{\rm init})_{S_{\rm init}^{\delta} } -
b^0 \|_1 =  \| (b^0 - \beta_{\rm init})_{S_{\rm init}^{\delta}}  \|_1 +
\| (b^0)_{S_0 \backslash S_{\rm init}^{\delta} }\|_1 . $$
Now,
$$ \| (b^0 - \beta_{\rm init})_{S_{\rm init}^{\delta}} \|_1 \le
\delta_1 $$
Moreover
$$ \| (b^0)_{S_0 \backslash S_{\rm init}^{\delta}} \|_1 \le
 \| (b^0 - \beta_{\rm init})_{S_0 \backslash S_{\rm init}^{\delta}} \|_1+
 \|( \beta_{\rm init})_{S_0 \backslash S_{\rm init}^{\delta}} \|_1 $$
 $$ \le  \| (b^0 - \beta_{\rm init})_{S_0 \backslash S_{\rm init}^{\delta}} \|_1+
\delta s_0  \le \delta_1 + \delta s_0 . $$
Hence
$$  \| ( \beta_{\rm init})_{S_{\rm init}^{\delta} }-
b^0 \|_1 \le 2 \delta_1 + \delta s_0 .$$
The $\ell_2$-error of the second result follows by the same arguments.

The first inequality  of the third result follows from the
definition of ${\rm f}_{S_{\rm init}^{\delta}}$ as projection, and the second follows
from the triangle inequality, where we invoke
that
$$| S_{\rm init}^{\delta} \backslash S_0 | \le {\delta_2^2 \over  \delta^2 }  $$
so that 
$$| S_{\rm init}^{\delta} | \le {\delta_2^2 \over  \delta^2 }+ s_0 , $$
and thus
$$\Lambda_{\rm max}^2 ( S_{\rm init}^{\delta} ) \le \biggl \lceil { \delta_2^2 \over \delta^2s_0 } +1 \biggr \rceil
 \Lambda_{\rm sparse}^2 (s_0) . $$
 The final result follows from
 $$\Lambda_{\rm min} ( S_{\rm init}^{\delta} \cup S_0) \ge
 \phi_{\rm sparse} ( S_0, | S_{\rm init}^{\delta} \backslash S_0 | + s_0 ) \ge
 \phi_{\rm sparse} ( S_0,2 s_0 ) . $$

\hfill $\sqcup \mkern -12mu \sqcap$

{\bf Proof of Theorem \ref{thres.theorem}.}

Inserting the bound $\delta_2 = O (\lambda_{\rm init} \sqrt {s_0} / \phi^2 (2, S_0 , 2 s_0) )$
(see Theorem \ref{init.theorem}),
and $\| {\rm f}_{S_0} - {\bf f}^0 \|= O (\lambda_{\rm init} \sqrt {s_0} / \phi^2 (2, S_0 ) )$,
we get for $\lambda_{\rm init} / \phi^2 ( 2 , S_0) = O ( \delta) $,  
$\delta \ge \delta_2 / \sqrt {s_0}$,
$$\| f_{S_{\rm init}^{\delta} } - {\bf f}^0 \|^2 =  \Lambda_{\rm sparse}^2 (s_0)
  \biggl [ { 1 \over \phi^4 ( 2, S_0 , 2s_0) }  +
 {\delta^2 \over \lambda_{\rm init}^2 } \biggr ] O (\lambda_{\rm init}^2 s_0)   , $$
 $$\| b^{S_{\rm init}^{\delta}} - b^0 \|_2  =
  {\Lambda_{\rm sparse} (s_0) 
 \over \phi_{\rm sparse}
 ( S_0 , 2s_0 ) } \times $$ $$\ \ \ \ \ \ \ \ \ \ \ 
  \biggl [ { 1 \over \phi^2 ( 2, S_0 , 2s_0) }  +
 {\delta \over \lambda_{\rm init} }  \biggr ] O( \lambda_{\rm init} \sqrt {s_0} ) , $$
 and
 $$| S_{\rm init}^{\delta} \backslash S_0 | = \biggl [
 { \lambda_{\rm init}^2 \over \delta^2 \phi^4 ( 2 , S_0 , 2s_0) } \biggr ] 
 O(s_0) . $$
 \hfill $\sqcup \mkern -12mu \sqcap$

\subsubsection{Proofs for Subsection \ref{adap.section}: the noiseless adaptive Lasso}

We use that when $\delta \ge \delta_2 /\sqrt {s_0}$, then
$S_{\rm init}^{\delta } \backslash S_0 \le s_0$. Application of Corollary
\ref{ordered.corollary}  then gives

\begin{corollary} \label{adap.corollary}
We have, for all $\delta \ge \delta_2 / \sqrt {s_0} $,  and all $\beta$
$$ \delta_{\rm adap}^2
 \le   2 \| f_{{\beta}_{S_{\rm init}^{\delta}}} - {\bf f}^0 \|^2 + 
{ 12  \lambda_{\rm init}^2 \lambda_{\rm adap}^2  s_0 \over
\delta^2 \phi_{\rm min}^2 (2, S_0 , 2 s_0) }        ,$$
and 
$$\| \beta_{\rm adap} - \beta_{S_{\rm init}^{\delta}}\|_1 \le 
{3 \delta \| f_{\beta_{S_{\rm init}^{\delta}}} - {\bf f}^0 \|^2 \over
 \lambda_{\rm init} \lambda_{\rm adap}} +
 {6 \lambda_{\rm init} \lambda_{\rm adap} s_0 \over
 \delta \phi_{\rm min}^2 ( 2 , S_0 , 2 s_0 ) } , $$
 and
$$\| \beta_{\rm adap} - \beta_{S_{\rm init}^{\delta}} \|_2 
\le { 6 \delta  \| f_{\beta_{S_{\rm init}^{\delta}} }  - {\bf f}^0 \|^2 
\over \sqrt {s_0} \lambda_{\rm init} \lambda_{\rm adap} }  + { 12 \lambda_{\rm init} \lambda_{\rm adap} \sqrt {s_0}  \over \delta \phi_{\rm min}^2 (2, S_0 ,
3 s_0) }  ,
 $$
 and, from Lemma \ref{betterthres},
 $$
 \| f_{(\beta_{\rm init})_{S_{\rm init}^{\delta}}} - {\bf f}^0 \|^2 \le 
2 \| {\rm f}_{S_0} - {\bf f}^0 \|^2 + 36 \Lambda_{\rm sparse}^2  (s_0)   
 \delta^2 s_0. $$
 Furthermore, from Lemma \ref{weightselect} ,
 $$|S_{\rm adap} \backslash S_0 |^2 \le
 \Lambda_{\rm max}^2  ( S_{\rm adap} \backslash S_0) 
 { \delta_{\rm adap}^2 \over
\lambda_{\rm adap}^2} {\delta_2^2 \over \lambda_{\rm init}^2 } . $$
  If $|S_{\rm adap} \backslash S_0| > s_0$, we have
$$|S_{\rm adap} \backslash S_0 | \le 8 \Lambda_{\rm sparse}^2 (s_0) { \delta_{\rm adap}^2 \over
\lambda_{\rm adap}^2 s_0} {\delta_2^2 \over \lambda_{\rm init}^2 } \wedge
 2 \Lambda_{\rm max} { \delta_{\rm adap} \over \lambda_{\rm adap}}
{\delta_2 \over \lambda_{\rm init}}.$$
 \end{corollary}
 
We note that in the above corollary, the use of the $\ell_2$-error 
 $\delta_2$ is rather crucial for the variable selection result: 
 with the weights $w_j = 1/ | \beta_{j,{ \rm init} } |$, we have
 $$\| ( 1/ w)_{S\backslash S_0} \|_2 =
 \| (\beta_{\rm init})_{S \backslash S_0} \|_2 \le \delta_2 . $$
 With alternative weights
 $w_j = 1/ \sqrt{ | \beta_{j, {\rm init} } | } $. 
 the theory can also be developed using only the $\ell_1$-error
 $\delta_1$.

A further observation is
that the above corollary is an obstructed oracle inequality, where the oracle is restricted 
to choose the index set as a thresholded set of the initial Lasso.
Concentrating on prediction error, it leads to defining the ``oracle" threshold as
 \begin{equation}\label{oraclethreshold}
\delta_0 := \arg \min_{\delta \ge \delta_2 / \sqrt {s_0} } 
\left \{  \| {\rm f}_{S_{\rm init}^{\delta}} - {\bf f}^0 \|^2 + 
{ 12 \lambda_{\rm init}^2   \lambda_{\rm adap}^2 s_0 \over
\delta^2 \phi_{\rm min}^2(2, S_0 , 2 s_0)}   
 \right \}  .
 \end{equation}
  This oracle
has active set $S_{\rm init}^{\delta_0}$, with size
$| S_{\rm init}^{\delta_0} | = O(s_0)$. Our following considerations however
will not be based on this optimal threshold, but rather on thresholds that
allow a comparison with the results for the thresholded initial Lasso.
This means that we might loose here some further favorable properties of the adaptive
Lasso.

 {\bf Proof of Theorem \ref{adap.theorem}.} 
 Corollary \ref{adap.corollary} combined with Lemma \ref{betterthres}
 gives that
for all $\delta \ge \delta_2 / \sqrt {s_0} $,
$$\delta_{\rm adap}^2 \le 
4  \| {\rm f}_{S_0} - {\bf f}^0 \|^2 + 72 \Lambda_{\rm sparse}^2 (s_0) \delta^2 s_0 
 +   { 12 \lambda_{\rm init}^2  \lambda_{\rm adap}^2 s_0  \over
\delta^2 \phi_{\rm min} ^2(2,S_0 , 2 s_0)   }    .$$
Using moreover that $ \| \beta_{\rm adap} - b^0 \|_q \le \| \beta_{\rm adap} - \beta_{S_{\rm init}^{\delta}}\|_q+ \| \beta_{S_{\rm init}^{\delta}} - b^0 \|_q$
and the bound of Lemma \ref{betterthres}, 
we get for $\delta \ge \delta_2 / \sqrt {s_0} $, 
$$\| \beta_{\rm adap} - b^0\|_1 \le { 3 \delta } s_0 +
{6 \delta \| {\rm f}_{S_0} -  {\bf f}^0 \|^2 \over
 \lambda_{\rm init} \lambda_{\rm adap}} 
 + {108 \Lambda_{\rm sparse}^2 (s_0) \delta^3 s_0  \over \lambda_{\rm init}
 \lambda_{\rm adap} } +
 {6 \lambda_{\rm init} \lambda_{\rm adap} s_0 \over
 \delta \phi_{\rm min}^2 ( 2 , S_0 , 2 s_0 ) } , $$
 and
$$\| \beta_{\rm adap} - b^0\|_2 
\le 3 \delta \sqrt {s_0} +{ 12 \delta  \| {\rm f}_{S_0}  - {\bf f}^0 \|^2 
\over \sqrt {s_0} \lambda_{\rm init} \lambda_{\rm adap} }  $$ $$+
{216 \Lambda_{\rm sparse}^2 (s_0) \delta^3 \sqrt {s_0}  \over \lambda_{\rm init}
 \lambda_{\rm adap} } + { 9 \lambda_{\rm init} \lambda_{\rm adap} \sqrt {s_0}  \over \delta \phi_{\rm min}^2 (2, S_0 ,
3 s_0) }  .
 $$
 Finally, again for $\delta \ge \delta_2 / \sqrt {s_0} $, 
 $$|S_{\rm adap} \backslash S_0 | \le $$ $${8 \Lambda_{\rm sparse}^2 (s_0)\delta_2^2 \over
\lambda_{\rm init}^2 \lambda_{\rm adap}^2  }  \left ( 
 { 4   \| {\rm f}_{S_0} - {\bf f}^0\|^2  \over  s_0 } +
{ 72  \Lambda_{\rm sparse}^2 (s_0) \delta^2  }   +  { 12   \lambda_{\rm init}^2 \lambda_{\rm adap}^2  \over
\delta^2   \phi_{\rm min}^2(2, S_0 , 2 s_0) }   \right )
 .  $$

 By Corollary \ref{init.corollary},
 $${\delta_2 \over \sqrt {s_0} } = O\biggl ( { \lambda_{\rm init} \over
 \phi^2 ( 2, S_0 , 2 s_0) } \biggr ) . $$
 Taking
 \begin{equation} \label{choicedelta}
  \delta^2 \asymp{  { \lambda_{\rm init} \lambda_{\rm adap} } 
 \over \phi_{\rm min} (2 , S_0 , 2 s_0)\Lambda_{\rm sparse} (s_0) } , 
 \end{equation}
 the requirement that $\delta \ge \delta_2 / \sqrt {s_0}$ is fulfilled if  take
 $$ \lambda_{\rm init}\biggl [  { \phi_{\rm min} ( 2 , S_0 , 2s_0) \Lambda_{\rm sparse}
 (s_0) \over  \phi^4 ( 2 , S_0 , 2s_0) } \biggr ] = O_{\rm suff} (\lambda_{\rm adap} ) ,$$
 that is, if Condition b holds.
We then obtain
 $$\delta_{\rm adap}^2 = \biggl [ { \Lambda_{\rm sparse} (s_0) \over
 \phi_{\rm min} ( 2 , S_0 , 2 s_0) }  \biggr ] O ( \lambda_{\rm init} 
 \lambda_{\rm adap} s_0 ) , $$
 $$ \|\beta_{\rm adap} - b^0 \|_1 =
 \biggl [ { \Lambda_{\rm sparse}^{1/2} (s_0) \over
 \phi_{\rm min}^{3/2} ( 2 , S_0 , 2 s_0) }  \biggr ] O (\sqrt { \lambda_{\rm init} 
 \lambda_{\rm adap} }s_0 ) , $$
  $$ \|\beta_{\rm adap} - b^0 \|_2 =
 \biggl [ { \Lambda_{\rm sparse}^{1/2}  (s_0) 
 \phi_{\rm min}^{1/2} (2 , S_0 ,2 s_0) \over
 \phi_{\rm min}^2 ( 2 , S_0 , 3 s_0) }  \biggr ] O (\sqrt { \lambda_{\rm init} 
 \lambda_{\rm adap}s_0 } ) , $$
 and
 $$| S_{\rm adap} \backslash S_0|  = 
{  \Lambda_{\rm sparse}^2 (s_0) \over
\phi^4 ( 2 , S_0 , 2 s_0) }  \biggl [  {\Lambda_{\rm sparse} (s_0) \over
\phi_{\rm min} ( 2 , S_0 , 2 s_0) } \biggr ] {\lambda_{\rm init} 
\over \lambda_{\rm adap} } O ( s_0) .  $$

\hfill $\sqcup \mkern -12mu \sqcap$

\subsubsection{Proofs for Subsection \ref{irrepresentable.section} on the weighted
irrepresentable condition}

{\bf Proof of Lemma \ref{irrepresentable}.}
This is the weighted variant of the first part of Lemma 6.2 in
\cite{vdG:2009}.
\hfill $\sqcup \mkern -12mu \sqcap$

{\bf Proof of Lemma \ref{boundirr.lemma}.}
We define, as in \cite{vdG:2009}, the {\it adaptive restricted regression}
$$\vartheta_{\rm adaptive} (S) :=
\max_{\beta \in {\cal R} (1,S)}
{ | ( f_{\beta_{S^c}} , f_{\beta_S}) | \over 
\| f_{\beta_S} \|^2 } . $$
Here, $(f, \tilde f )$ denotes the inner product between $f $ and $\tilde f$ as
elements of $L_2 (Q)$.

We will show that
\begin{equation}\label{weshow}
\sup_{\| \tau_{S} \|_{\infty} \le 1 } 
 \| W_{S^c}^{-1} \Sigma_{2,1} (S) \Sigma_{1,1}^{-1} (S)
 W_{S} \tau_{S} \|_{\infty} \le { \| w_{S} \|_2 \over \sqrt {|S|} 
 w_{S^c}^{\rm min} } \vartheta_{\rm adaptive} (S)  .
 \end{equation}
 
 It is moreover not difficult to see that $\vartheta_{\rm adaptive}(S) \le \sqrt {|S|} / 
\Lambda_{\rm min} ( S) $, so then the proof of Lemma \ref{boundirr.lemma}
is done. 

 To derive (\ref{weshow}), we first note that
  $$\| W_{S^c}^{-1}\Sigma_{2,1} (S) \Sigma_{1,1}^{-1} (S)
 W_{S} \tau_{S} \|_{\infty} \le   \| \Sigma_{2,1} (S) \Sigma_{1,1}^{-1} (S)
 W_{S} \tau_{S} \|_{\infty} / w_{S^c}^{\rm min} .$$ 
 
 Define
 $$\beta_{S} := \Sigma_{1,1}^{-1} (S) W_{S} \tau_{S} .$$
  Then
 $$ \|W_{S^c}^{-1}  \Sigma_{2,1} (S) \Sigma_{1,1}^{-1} (S)
 W_{S} \tau_{S} \|_{\infty}=
 \sup_{\| \gamma_{S^c} \|_1 \le 1 }
|  \gamma_{S^c}^T W_{S^c}^{-1}  \Sigma_{2,1} (S) \Sigma_{1,1}^{-1} (S)
 W_{S} \tau_{S} |$$
 $$ = \sup_{\| W_{S^c} \beta_{S^c} \|_1 \le 1}
 | \beta_{S}^T \Sigma_{2,1} (S) \beta_{S} | =
 \sup_{\| W_{S^c} \beta_{S^c} \|_1 \le 1} |( f_{\beta_{S^c}} ,
 f_{\beta_{S}} ) | $$
 $$ \le \sup_{\|\beta_{S^c} \|_1 \le 1/w_{S^c}^{\rm min}} |( f_{\beta_{S^c}} ,
 f_{\beta_{S}} ) | $$
 $$ = \sup_{\|  \beta_{S^c} \|_1 \le \| w_{S} \|_2 \| \beta_{S} \|_2 /
 w_{S^c}^{\rm min} }
{  |( f_{\beta_{S^c}} ,
 f_{\beta_{S}} ) | \over  \| w_{S} \|_2 \| \beta_{S} \|_2 } $$
 $$= \sup_{\| \beta_{S^c} \|_1 \le  \| w_{S} \|_2 \| \beta_{S} \|_2  / w_{S^c}^{\rm min} }
{  |( f_{\beta_{S^c}} ,
 f_{\beta_{S}} ) | \over  \| f_{\beta_{S}} \|^2 } 
{   \| f_{\beta_{S}} \|^2  \over   \| w_{S}   \|_2 \| \beta_{S} \|_2  } . 
 $$
 But
 $$
{  \| f_{\beta_{S}} \|^2  \over   \| w_{S}   \|_2 \| \beta_{S} \|_2  } =
{ \tau_{S}^T W_{S} \Sigma_{1,1}^{-1} (S) W_{S}
\tau_{S} \over
\sqrt {\tau_{S}^T  W_{S}^2 \tau_{S} }\sqrt { \tau_{S} W_{S} \Sigma_{1,1}^{-2}
(S) 
W_{S} \tau_{S} } } { \| W_{S} \tau_{S} \|_2 \over
\| w_{S} \|_2 } \le 1 . $$
We conclude that
$$ \|W_{S^c}^{-1}  \Sigma_{2,1} (S) \Sigma_{1,1}^{-1} (S)
 W_{S} \tau_{S} \|_{\infty}\le 
 \sup_{\| \beta_{S^c} \|_1 \le  \| w_{S} \|_2\| \beta_{S} \|_2  / w_{S^c}^{\rm min} }
{  |( f_{\beta_{S^c}} ,
 f_{\beta_{S}} ) | \over  \| f_{\beta_{S}} \|^2 } $$ $$=
  { \| w_{S} \|_2 \over \sqrt {|S|} 
 w_{S^c}^{\rm min} } \vartheta_{\rm adaptive} (S).$$
 \hfill $\sqcup \mkern -12mu \sqcap$

\subsection{Proofs for Section \ref{noise.section}: the noisy case
}
Theorem \ref{weightnoise} gives bounds for prediction error, estimation
error and the number of false positives of the noisy weighted
Lasso. 

{\bf Proof of Theorem \ref{weightnoise}.} 
We can derive the prediction and estimation results in the same way as in Theorem \ref{weight},
adding now the noise term: 
$$\|  \hat f_{\rm weight} - {\bf f}^0 \|_n^2 +\lambda_{\rm init}
\lambda_{\rm weight}  \sum_{j=1}^p w_j | \hat \beta_{j, {\rm weight}}  | $$ $$\le 
2 ( \epsilon , \hat f_{\rm weight} - f_{\beta_S} )_n + \| f_{\beta_S} - {\bf f}^0\|_n^2 +
\lambda_{\rm init}
\lambda_{\rm weight} \sum_{j \in S } w_j |\beta_j | $$
$$ \le \lambda_{ \rm init} \|\hat  \beta_{\rm weight}  - 
\beta_S \|_1 /2 + \| f_{\beta_S} - {\bf f}^0\|_n^2 +
\lambda_{\rm init}
\lambda_{\rm weight} \sum_{j \in S } w_j |\beta_j | $$
and hence, using $\lambda_{\rm weight} w_{S^c}^{\min} \ge 1 $,
$$\|  \hat f_{\rm weight} - {\bf f}^0 \|_n^2 +\lambda_{\rm init}
\lambda_{\rm weight} w_{S^c}^{\min} 
\| \hat \beta_{S^c} \|_1  /2 $$ $$\le 
  \| f_{\beta_S} - {\bf f}^0\|_n^2 +\biggl [ \lambda_{\rm int} / 2 + \lambda_{\rm init}
  \lambda_{\rm weight} \| w_S \|_2 /\sqrt {s} \biggr ] \sqrt {s} \|
  \hat \beta_{\rm weight} - \beta_S \|_2  
.$$

As  $\lambda_{\rm weight} \| w_S \|_2/\sqrt {s} \ge 1 $
it gives
$$\| \hat  f_{\rm weight} - {\bf f}^0 \|_n^2 +\lambda_{\rm init}
\lambda_{\rm weight} w_{S^c}^{\rm min} 
\| \hat \beta_{S^c} \|_1  /2 $$ $$\le 
  \| f_{\beta_S} - {\bf f}^0\|_n^2 +3 \lambda_{\rm int}  \lambda_{\rm weight}
  \| w_S \|_2    \|
  \hat \beta_{\rm weight} - \beta_S \|_2 /2.$$
  Now insert $w_{S^c}^{\rm min} \ge M/L$,
  $1 \le \lambda_{\rm weight } M$ and $\| w_S \|_2 /\sqrt s \le M$:
  $$\|  \hat f_{\rm weight} - {\bf f}^0 \|_n^2 +\lambda_{\rm init}
\lambda_{\rm weight} M 
\|\hat  \beta_{S^c} \|_1  /(2L) $$ $$\le 
  \| f_{\beta_S} - {\bf f}^0\|_n^2 +3 \lambda_{\rm int}  \lambda_{\rm weight}
M \sqrt {s} 
  \| \hat \beta_{\rm weight} - \beta_S \|_2 /2.$$
  The rest of the proof for the prediction and estimation error can therefore carried out in the same way is
the proof of Theorem \ref{weight}. 

As for variable selection, we use
as in Lemma \ref{weightselect} the weighted KKT conditions: for all $j$
$$2 ( \psi_j , \hat f_{\rm weight} - {\bf f}^0 )_n  -2 (\psi_j , \epsilon)_n 
=- \lambda_{\rm init} \lambda_{\rm weight} w_j \hat \tau_{j, {\rm weight}} ,$$
where $\| \hat \tau_{\rm weight} \|_{\infty} \le 1 $ and
$\hat \tau_{j, {\rm weight}} {\rm l} \{ \hat \beta_{j, {\rm weight}} \not= 0 \} =
{\rm sign} ( \hat \beta_{j,{\rm weight}} ) $.
Invoking $\lambda_{\rm weight} w_{S^c}^{\rm min} 
\ge 1$, we know that for all $j \in S^c$, $\lambda_{\rm weight} w_j \ge
1$. Moreover, $2 |(\epsilon , \psi_j )_n \le \lambda_{\rm init}/2$
by the definition of ${\cal T}$. Therefore, 
$$\sum_{j \in \hat S_{\rm weight} \cap S^c  \backslash S_0 }
2 | ( \psi_j , \hat f_{\rm weight}  - {\bf f}^0 )_n |^2 \ge \lambda_{\rm init}^2 
\lambda_{\rm weight}^2\| w_{\hat S_{\rm weight} \cap S^c \backslash S_0 } \|_2^2 /4
. $$
One can now proceed as in Lemma \ref{weightselect}.

\hfill $\sqcup \mkern -12mu \sqcap$

\subsubsection
{ Proof of Lemma \ref{refineselect.lemma} with the more involved
conditions} 
To prove this lemma, we actually need some results
in from Section \ref{preview.section} and an intermediate result in their
proof. One may skip the present proof at first reading and first consult
the next subsection
(Subsection \ref{preview.proofs}).

The bound for the number of false positives of the initial lasso follows from
the inequality
$$| \hat S_{\rm init} \backslash S_0 | \le 
{ \Lambda_{\rm max}^2  (  \hat S_{\rm init} \backslash S_0)  \over 
\phi^2 ( 6 , S_0 ) } O(s_0)  .$$
This follows from Theorem \ref{weightnoise}, and from inserting the
bound of Theorem \ref{noisyinit.theorem} for $\hat \delta_{\rm init}$.
One can then proceed by applying the inequality
\begin{equation} \label{s*}
\Lambda_{\rm max}^2  (  \hat S_{\rm init} \backslash S_0)  \le
\biggl ( { | \hat S_{\rm init} \backslash S_0 | \over s_* } + 1 \biggr )
\Lambda_{\rm sparse}^2 ( s_*) . 
\end{equation}
The result for the adaptive Lasso can be derived from
$$| \hat S_{\rm adap} \backslash S_0 |^2 \le 
{ \Lambda_{\rm max}^2  (  \hat S_{\rm adap} \backslash S_0)  \over 
\phi^4 ( 6 , S_0, 2 s_0 ) } \biggl [ { \Lambda_{\rm sparse} (s_0) \over
\phi_{\rm min} (6, S_0 , 2 s_0) } \biggr ] {\lambda_{\rm init} \over
\lambda_{\rm adap} } O(s_0)  .$$
This follows from  (\ref{adap*}) (which can be found
at the end of the proof of Theorem \ref{noisyadap.theorem}), invoking Condition B,
and applying the bound of Theorem \ref{noisyadap.theorem} for
$\hat \delta_{\rm adap}$, and the bound of Theorem
\ref{noisyinit.theorem} for $\hat \delta_2 $. 
Insert again (\ref{s*}) to complete the proof. 
\hfill $\sqcup \mkern -12mu \sqcap$

\subsection{Proofs for Section \ref{preview.section}} \label{preview.proofs}

\subsubsection{Proof of  the probability inequality of Lemma \ref{noise.lemma}}
This follows easily from
the probability bound $\PP (|Z| \ge \sqrt {2t}) \le 2 \exp[-t]$ for a standard 
normal random variable $Z$.
\hfill $\sqcup \mkern -12mu \sqcap$

\subsubsection{Proof of 
Theorem \ref{noisyinit.theorem}: the noisy initial Lasso}

Theorem \ref{noisyinit.theorem} is a simplified formulation of Corollary
\ref{noisyinit.corollary} below. This corollary follows from
Theorem \ref{weightnoise} by taking $L=1$ and $S= S_0$.

\begin{corollary}\label{noisyinit.corollary}
Let
$$ \delta_{\rm oracle}^2 := 
 \| {\rm f}_{S_0} - {\bf f}^0 \|_n^2 + { 7 \lambda_{\rm init}^2 |S_0| \over \phi^2 (6, S_0, 2 s_0)  } . $$
 Take $\lambda_{\rm init} \ge 2 \lambda_{\rm noise}$.
We have on ${\cal T}$,
$$ \hat \delta_{\rm init}^2 \le
2 \delta_{\rm oracle}^2.
$$
Moreover, on ${\cal T}$,
$$ \hat \delta_1 \le 
5 \delta_{\rm oracle}^2 / \lambda_{\rm init} , $$
and
$$\hat  \delta_2 \le  10 \delta_{\rm oracle}^2/ ( \lambda_{\rm init} \sqrt {s_0})  . $$
Also, on ${\cal T}$, 
$$| \hat S_{\rm init} \backslash S_0 | \le 16 \Lambda_{\rm max}^2 (\hat S_{\rm init}
\backslash S_0) 
 { \hat  \delta_{\rm init}^2 \over \lambda_{\rm init}^2} .$$
\end{corollary}

\subsubsection{Proof of Theorem \ref{noisythres.theorem}: the noisy thresholded Lasso}

The least squares estimator $\hat {\rm f}_{\hat S_{\rm init}^{\delta}}$ using only variables in $\hat S_{\rm init}^{\delta}$ (i.e., the projection of ${\bf Y}= {\bf f}^0 + \epsilon$ on the
linear space
spanned by $\{ \psi_j \}_{j \in \hat S_{\rm init}^{\delta}} $)
has similar prediction properties as
${\rm f}_{\hat S_{\rm init}^{\delta} }$ (the projection of ${\bf f}^0$ on the same
linear space).
This is because, as is shown
in the next lemma, their difference is small.

\begin{lemma} \label{leastsquares}
Let $\delta \ge \hat \delta_2 / \sqrt {s_0} $. Then on ${\cal T}$,
$$\| \hat {\rm f}_{\hat S_{\rm init}^{\delta}} - {\rm f}_{\hat S_{\rm init}^{\delta } } \|_n^2 \le 
 {  \lambda_{\rm init}^2  s_0 \over 2 \phi_{\rm sparse}^2 (S_0 , 2 s_0) } . $$
\end{lemma}

{\bf Proof of Lemma \ref{leastsquares}.} 
This follows from
$$\| \hat {\rm f}_{\hat S_{\rm init}^{\delta}} - {\rm f}_{\hat S_{\rm init}^{\delta } } \|_n^2 
\le 2 ( \epsilon , \hat {\rm f}_{\hat S_{\rm init}^{\delta}} - {\rm f}_{\hat S_{\rm init}^{\delta } })_n , $$
and
$$2 ( \epsilon , \hat {\rm f}_{\hat S_{\rm init}^{\delta}} - {\rm f}_{\hat S_{\rm init}^{\delta } })_n  \le \lambda_{\rm init} 
\| \hat b^{\hat S_{\rm init}^{\delta}} - b^{\hat S_{\rm init}^{\delta}} \|_1 /2$$
$$ \le \lambda_{\rm init}\sqrt {2 s_0} \| \hat b^{\hat S_{\rm init}^{\delta}} - b^{\hat S_{\rm init}^{\delta}} \|_2 /2
 \le  \lambda_{\rm init} \sqrt {2 s_0} \| \hat {\rm f}_{\hat S_{\rm init}^{\delta}} - {\rm f}_{\hat S_{\rm init}^{\delta } } \|_n/(2 \phi_{\rm sparse} (S_0 , 2 s_0))  . $$
\hfill $\sqcup \mkern -12mu \sqcap$

{\bf Proof of Theorem \ref{noisythres.theorem}}
The bound for 
$\| (\hat \beta_{\rm init}){\hat S_{\rm init}^{\delta}} -
b^0 \|_2 \le 2 \hat \delta_2 + \delta \sqrt {s_0}$ can be derived in the same way as
in Lemma \ref{betterthres}. The same is true for the bound
$$ 
\| {\rm f}_{\hat S_{\rm init}^{\delta} }  - {\bf f}^0 \|_n \le
\| f_{(\hat \beta_{\rm init})_{\hat S_{\rm init}^{\delta}}} - {\bf f}^0 \|_n $$ $$ \le 
 \| {\rm f}_{S_0} - {\bf f}^0 \|_n + \sqrt {\biggl  \lceil {\hat \delta_2^2 \over  \delta^2 s_0}
 + 1 \biggr \rceil }   \Lambda_{\rm sparse}  (s_0) ( 2\hat  \delta_2 + \delta  \sqrt{ s_0}).
  $$
  Assumption A together with Lemma \ref{leastsquares} complete
  the proof for the bounds for prediction and estimation error,
  with the $\ell_1$-bound being a simple consequence of the $\ell_2$-bound.
Also, the variable selection result follows from
$$| \hat S_{\rm init}^{\delta } \backslash S_0 | \le
{\hat \delta_2^2 \over \delta^2 } , $$
and Assumption A.
\hfill $\sqcup \mkern -12mu \sqcap$

\subsubsection{Proof of Theorem \ref{noisyadap.theorem}: the noisy adaptive Lasso}

We first apply Theorem \ref{weightnoise} to the adaptive Lasso.

\begin{corollary} \label{noisyadap.corollary}
Suppose we are on ${\cal T}$. Take 
$\lambda_{\rm adap} \ge \delta \ge \hat \delta_2 / \sqrt {s_0} $. 

We have, for all $\delta \ge \hat \delta_2 / \sqrt {s_0} $,  and all $\beta$
$$\hat \delta_{\rm adap}^2
 \le   2 \| f_{{\beta}_{\hat S_{\rm init}^{\delta}}} - {\bf f}^0 \|_n^2 + 
{ 28  \lambda_{\rm init}^2 \lambda_{\rm adap}^2  s_0 \over
\delta^2 \phi_{\rm min}^2 (6, S_0 , 2 s_0) }        ,$$
and 
$$\| \hat \beta_{\rm adap} - \beta_{\hat S_{\rm init}^{\delta}}\|_1 \le 
{5 \delta \| f_{\beta_{\hat S_{\rm init}^{\delta}}} - {\bf f}^0 \|_n^2 \over
 \lambda_{\rm init} \lambda_{\rm adap}} +
 {14 \lambda_{\rm init} \lambda_{\rm adap} s_0 \over
 \delta \phi_{\rm min}^2 ( 6 , S_0 , 2 s_0 ) } , $$
 and
$$\| \hat \beta_{\rm adap} - \beta_{\hat S_{\rm init}^{\delta}} \|_2 
\le { 10 \delta  \| f_{\beta_{\hat S_{\rm init}^{\delta}} }  - {\bf f}^0 \|_n^2 
\over \sqrt {s_0} \lambda_{\rm init} \lambda_{\rm adap} }  + { 42 \lambda_{\rm init} \lambda_{\rm adap} \sqrt {s_0}  \over \delta \phi_{\rm min}^2 (6, S_0 ,
3 s_0) }  .
 $$
Moreover
$$|( \hat S_{\rm adap}\cap (\hat S_{\rm init}^{\delta})^c )  \backslash S_0 | \le  s_0 + 32
\Lambda_{\rm sparse} (s_0)  { \hat \delta_{\rm adap}^2 \over
\lambda_{\rm adap}^2s_0 } {\hat \delta_2^2 \over \lambda_{\rm init}^2 } \wedge
 4 \Lambda_{\rm max} { \hat \delta_{\rm adap} \over \lambda_{\rm adap}}
{\hat \delta_2 \over \lambda_{\rm init}}. $$
\end{corollary}

{\bf Proof of Theorem \ref{noisyadap.theorem}.}

By the same arguments as used in Lemma \ref{betterthres},
 for $\delta \ge \hat \delta_2 /\sqrt {s_0}$,
 $$
 \| f_{(\hat \beta_{\rm init})_{\hat S_{\rm init}^{\delta}}} - {\bf f}^0 \|_n \le 
 \| {\rm f}_{S_0} - {\bf f}^0 \|_n^2 + 3 \sqrt 2  \Lambda_{\rm sparse}^2  (s_0)   
 \delta^2 s_0, $$
and
$\| (\hat \beta_{\rm init})_{\hat S_{\rm init}^{\delta}} -
b^0 \|_2 \le 3 \delta \sqrt {s_0}$. 
The prediction and estimation results now follow from Corollary
\ref{noisyadap.corollary} combined with Condition B. 

We apply Corollary \ref{noisyadap.corollary} with 
\begin{equation}\label{delta}
\delta^2 = { \lambda_{\rm init} \lambda_{\rm adap}
\over \phi_{\rm min} ( 6, S_0 , 2 s_0) \Lambda_{\rm sparse} } . 
\end{equation}
Condition B requires that
$$ \biggl [ { \Lambda_{\rm sparse}
(s_0) 
 \over \phi_{\min}^3 ( 6 , S_0 , 2 s_0) }\bigg ] \lambda_{\rm init}= O_{\rm suff} (\lambda_{\rm adap}) . $$
This ensures that $\delta \ge \hat \delta_2 / \sqrt {s_0}$
on the set ${\cal T}$. Moreover, equation (\ref{delta}) gives that
$\lambda_{\rm adap}\ge \delta$ as soon as
$$\lambda_{\rm adap} \ge \biggl [ {1 \over \phi_{\rm min} ( 6, S_0 , 2 s_0)
\Lambda_{\rm sparse}(s_0) } \biggr ] \lambda_{\rm init} , $$
which  is also ensured by Condition B.

The variable selection result follows from: for 
$\delta \ge \hat \delta_2 / \sqrt {s_0}$,
\begin{equation} \label{adap*}
| \hat S_{\rm adap} \backslash S_0| \le
| ( \hat S_{\rm adap} \cap ( \hat S_{\rm init}^{\delta})^c   \backslash S_0| +
| \hat S_{\rm init}^{\delta } \backslash S_0 |\le
| ( \hat S_{\rm adap} \cap ( \hat S_{\rm init}^{\delta})^c   \backslash S_0|  + s_0 .
\end{equation}
\hfill $\sqcup \mkern -12mu \sqcap$

\subsubsection{Proof of Lemma
\ref{min.lemma}, where coefficients are assumed to be large}
On ${\cal T}$, 
for $j \in S_0$, $| \hat \beta_{j, {\rm init}} | >\hat  \delta_{\infty}$,
and $| \hat \beta_{j, {\rm init}} | >| b_j^0|/2 $,
since $|b_j^0 | > 2 \hat \delta_{\infty}$.
Moreover, for $ j \in S_0^c$, $| \hat \beta_{j, {\rm init}} | \le \hat \delta_{\infty}$.
Let
$$M^2={ 4 \over s_0}
 \sum_{j \in S_0 }
{1 \over | b_j^0|^2 } . $$
So
$$\| w_{S_0} \|_2^2 / s_0 \le M^2  . $$
Note that $M \le 1/ \hat \delta_{\infty}$. 
Since $w_{S_0^c}^{\min} \ge 1/\hat \delta_{\infty}$,
the condition $\lambda_{\rm adap} M \ge 1 $ implies
$\lambda_{\rm adap} w_{S_0^c}^{\rm min} \ge 1$.

Apply Theorem \ref{weightnoise} to the adaptive Lasso with $S=S_0$,
and $\beta = b^0$:
$$\hat \delta_{\rm adap}^2
 \le   2 \| {\rm f}_{S_0} - {\bf f}^0 \|_n^2 + 
{ 14  \lambda_{\rm init}^2 \lambda_{\rm adap}^2  M^2s_0 \over
 \phi^2 (6, S_0 ) }   =
O \biggl (  { \lambda_{\rm init}^2 \lambda_{\rm adap}^2  M^2s_0 \over
 \phi^2 (6, S_0 ) }\biggr ) 
      ,$$
and 
$$\| \hat \beta_{\rm adap} -b^0\|_1 \le 
{5  \| {\rm f}_{S_0} - {\bf f}^0 \|_n^2 \over
 \lambda_{\rm init} \lambda_{\rm adap}M } +
 {7 \lambda_{\rm init} \lambda_{\rm adap} M s_0 \over
  \phi^2 ( 6 , S_0 , ) } = O \biggl (
   { \lambda_{\rm init} \lambda_{\rm adap} M s_0 \over
  \phi^2 ( 6 , S_0  ) } \biggr )  ,$$
 and
$$\| \hat \beta_{\rm adap} -b^0 \|_2 
\le { 10   \|{\rm f}_{S_0}   - {\bf f}^0 \|_n^2 
\over M \sqrt {s_0} \lambda_{\rm init} \lambda_{\rm adap} }  + { 28 \lambda_{\rm init} \lambda_{\rm adap} M  \sqrt {s_0}  \over  \phi^2 (6, S_0 ,
2 s_0) }  = O \biggl ( {  \lambda_{\rm init} \lambda_{\rm adap} M  \sqrt {s_0}  \over  \phi^2 (6, S_0 ,
2 s_0) } \biggr ) .
 $$

Also,  when $|\hat S_{\rm adap}  \backslash S_0 |> s_0$,
it holds that 
$$| \hat S_{\rm adap}   \backslash S_0 | \le 
32 \Lambda_{\rm sparse}^2 (s_0)  {\| \hat f_{\rm adap} - {\bf f}^0 \|_n^2   \over \lambda_{\rm adap}^2 s_0 }
{  \|( 1/w)_{\hat S_{\rm adap} \backslash S_0}\|_2^2 \over \lambda_{\rm init}^2}
$$  
$$ \le 32 \Lambda_{\rm sparse}^2 (s_0)  {\| \hat f_{\rm adap} - {\bf f}^0 \|_n^2   \over \lambda_{\rm adap}^2 s_0 }
{  \hat \delta_2^2  \over \lambda_{\rm init}^2} $$ $$=
 \Lambda_{\rm sparse}^2 (s_0) O \biggl (   { \lambda_{\rm init}^2   M^2 s_0  \over
 \phi^2 (6, S_0 )\phi^4 ( 6 , S_0 , 2 s_0)  }\biggr ) .
$$
\hfill $\sqcup \mkern -12mu \sqcap$

\bibliographystyle{plainnat}
\bibliography{reference}

\end{document}